# Honest variable selection in linear and logistic regression models via $\ell_1$ and $\ell_1 + \ell_2$ penalization

## Florentina Bunea[*]

*Department of Statistics, Florida State University, Tallahassee, Florida,*
*e-mail:* flori@stat.fsu.edu

**Abstract:** This paper investigates correct variable selection in finite samples via $\ell_1$ and $\ell_1 + \ell_2$ type penalization schemes. The asymptotic consistency of variable selection immediately follows from this analysis. We focus on logistic and linear regression models. The following questions are central to our paper: given a level of confidence $1 - \delta$, under which assumptions on the design matrix, for which strength of the signal and for what values of the tuning parameters can we identify the true model at the given level of confidence? Formally, if $\widehat{I}$ is an estimate of the true variable set $I^*$, we study conditions under which $\mathbb{P}(\widehat{I} = I^*) \geq 1 - \delta$, for a given sample size $n$, number of parameters $M$ and confidence $1 - \delta$. We show that in identifiable models, both methods can recover coefficients of size $\frac{1}{\sqrt{n}}$, up to small multiplicative constants and logarithmic factors in $M$ and $\frac{1}{\delta}$. The advantage of the $\ell_1 + \ell_2$ penalization over the $\ell_1$ is minor for the variable selection problem, for the models we consider here. Whereas the former estimates are unique, and become more stable for highly correlated data matrices as one increases the tuning parameter of the $\ell_2$ part, too large an increase in this parameter value may preclude variable selection.

Received August 2008.

## 1. Introduction

The literature on various theoretical aspects of $\ell_1$ empirical risk minimization has enjoyed substantial growth over the last few years, partly as a necessity to complement the flourishing field of convex optimization. The main attraction, from both theoretical and computational perspectives, is the proved ability of such methods to recover sparse approximations of the true underlying model when the number of parameters is large relative to the sample size. The principal theoretical topics of interest are therefore focused on optimality properties that involve the notion of sparsity. Whereas the theoretical properties of the $\ell_1 + \ell_2$ penalized estimates, sometimes referred to as elastic net estimates, a phrase introduced by [23] in linear models, have not been investigated for the models we consider, the properties of the $\ell_1$ penalized estimates, typically referred to as the Lasso-type estimates, have received considerable attention. The topics studied range from finite sample results concerning sparsity oracle inequalities for the risk of the estimators, in regression and classification, e.g., [4, 5, 19, 26, 20, 2, 11] to the asymptotic behavior of the estimates, including the consistency of subset selection, e.g. [9, 10, 13, 22, 21, 25, 6, 3, 17, 12, 14].

[*]Research partially supported by NSF grant DMS 0706829 and the Isaac Newton Institute, Cambridge, UK.

This work is motivated by the emergence of a large number of variations and improvements of the $\ell_1$ penalization schemes in regression and classification. To appreciate the need for such variations it is important therefore to investigate the limitations of the original method. When the number of variables $M$ is large relative to $n$, an asymptotic analysis of the variable selection problem may obscure issues that arise in finite samples. In this paper we investigate the finite sample accuracy of variable selection via the $\ell_1$ and the closely related $\ell_1 + \ell_2$ penalization schemes in regression models. We also discuss asymptotic alternatives and asymptotic consequences of our results. Our goal is to review existing results, and to offer a self-contained, back-to-back analysis of these important models and respective penalization schemes.

Formally, let $(\mathbf{X}_i, Y_i)$, $1 \leq i \leq n$, be i.i.d. pairs distributed as $(\mathbf{X}, Y)$ with probability measure $\mathbb{P}$, where $Y \in \{0, 1\}$ or $Y \in \mathbb{R}$ and $\mathbf{X} = (X_1, \ldots, X_M) \in \mathbb{R}^M$. We assume that $\mathbb{E}(Y|\mathbf{X} = x) = g(\sum_{j \in I^*} \beta_j^* x_j)$, where $I^* \subseteq \{1, \ldots, M\}$ is an unknown subset and $g$ is a known link function. In our analysis, $M$ is allowed to depend and be larger than the sample size $n$, and the size of $I^*$ may depend on $n$. The goal of this paper is to provide an understanding of the merits and possible limitations of variable selection via these two penalization schemes when used to answer the following central questions: given a level of confidence $1 - \delta$, given the number of variables $M$ and the sample size $n$, under which assumptions on the design matrix, for which strength of the signal and for what values of the tuning parameters do we identify the true model at the given level of confidence? Formally, if $\widehat{I}$ is an estimate of $I^*$, we study conditions under which $\mathbb{P}(\widehat{I} = I^*) \geq 1 - \delta$.

We will focus on variable selection in logistic regression, corresponding to the link function $g(z) = e^z/(1 + e^z)$, and also present a full analysis of the problem for linear models, corresponding to $g(z) = z$, to facilitate the comparison of the results. We will conduct separate analyses of the corresponding estimates, as different arguments are needed for models with possibly unbounded response, such as the linear model.

We denote by $\beta^*$ the vector in $\mathbb{R}^M$ with components $\beta_j^*$ for $j \in I^*$ and zero otherwise. We begin our analysis in Section 2 by establishing upper bounds on the $\ell_1$ distance between the Lasso and elastic net estimators, respectively, and the parameter $\beta^*$. These results are connected with the sparsity oracle inequalities recently obtained for the Lasso estimators in [4] and [2], in linear regression models, and [19], in generalized linear regression models. The focus in these works is on the predictive performance of the estimators, rather than on the accuracy of variable selection, as considered here. For us, these results are an intermediate, albeit essential, step in discussing the conditions under which an estimate $\widehat{I}$ of the set $I^*$ satisfies $\mathbb{P}(\widehat{I} = I^*) \geq 1 - \delta$. It is intuitively clear that if the estimates $\widehat{\beta}$ are too far from $\beta^*$, we cannot hope to recover the true coefficient set $I^*$ with high probability. It is interesting to note, however, that under some conditions on the design matrix, we can still estimate the true subset correctly even if the distance between $\beta^*$ and the estimates is not close to zero, but can still be controlled as in Section 2. Although this may appear surprising, it is this phenomenon that sets the variable selection problem apart from the problem of estimating well $\beta^*$ itself: here we aim at identifying a non-zero coefficient. Even if the estimate of this coefficient is relatively far from the real value, it only matters whether it is different than zero, not whether it is very close to the truth.

The rest of paper is organized as follows. In Section 2.1 we re-visit the conditions on the design matrix under which sparsity oracle inequalities for the Lasso estimates have been previously established and discuss weaker conditions. In Sections 2.2 and 2.3 we show that these results continue to hold under the weaker conditions. If one considers a slight modification of the $\ell_1$ penalty that consists in the addition of a properly scaled $\ell_2$ term, one can further weaken the requirements on the design matrix, while maintaining the sparsity of the resulting estimator. This motivates the study the $\ell_1 + \ell_2$ estimates, which have not been, to the best of our knowledge, investigated theoretically from this perspective in these models. Section 2.3 also contains an alternative asymptotic analysis of the $\ell_1$ norm of the difference $\widetilde{\beta} - \beta^*$ for estimates in logistic regression, motivated by the presence of possibly large constants in the finite sample oracle bounds. Under weak conditions on a weighted version of the design matrix, we obtain improved oracle bounds, that hold with probability converging to one.

In Section 3, which is central to our paper, we discuss in detail when the Lasso and the elastic net methods can provide accurate variable selection, in linear and logistic regression models. We show that obtaining results of the type $\mathbb{P}(I^* = \widehat{I}) \geq 1 - \delta$ depends crucially on a combination of conditions on the design matrix and the signal strength. This analysis complements the existing asymptotic results for Lasso estimates in linear regression models, and shows that similar phenomena occur in generalized linear models, for which the variable selection problem has not been investigated from this perspective; we refer to the the very recent work in [17] for related results in binary graphical models. Moreover, we provide the parallel study of the elastic estimates, and investigate to which extent they can be used for variable selection. We note that in a non-asymptotic framework, the study of $\mathbb{P}(I^* = \widehat{I})$ is well posed only if $\widehat{I}$ is unique. Since the elastic net estimates of $\beta^*$ are unique, as shown in Appendix B, so is the corresponding $\widehat{I}$. Recall that, in contrast, the Lasso-type estimators of $\beta^*$ may not be unique. However, in that case, the problem studied here is still well posed: even when the Lasso estimates of $\beta^*$ are not unique, the corresponding $\widehat{I}$ is. This property has been used implicitly in [15], and then in [13], for linear models, without an explicit proof, and not investigated outside linear models. For completeness, we present a proof of this result in Appendix B.

In the Conclusions section we summarize our findings and discuss the relative merits of the Lasso and elastic net estimates. The proofs of our main results are in Appendix A. Additional technical results are collected in Appendix B.

### 1.1. Notation

In the following sections we will denote by $\widehat{\beta}$ the penalized least squares estimates, for both the $\ell_1$ and $\ell_1 + \ell_2$ penalties and, similarly, by $\widetilde{\beta}$ the penalized logistic regression estimates, for either penalty. The estimates are of course different, but we opted for the same notation to keep the exposition simple. It will always be clear from the context to which combination model/penalty they correspond to. In the same way, $\widehat{I}$ will always denote the set of selected variables, and $I^*$ will denote the set of truly associated variables. We denote by $k^*$ the cardinality of $I^*$. For simplicity, we assume that the observations on the $X$ variables are normalized and centered, that is $\frac{1}{n}\sum_{i=1}^{n} X_{ij}^2 = 1$ and $\frac{1}{n}\sum_{i=1}^{n} X_{ij} = 0$, for all $j$. This is in no way crucial, but it allows for cleaner results and easier interpretation of the assumptions. We will also assume that for all $i$ and $j$ the variables $X_{ij}$ are bounded by a common constant $L > 0$, with probability 1. For any vector in $a \in \mathbb{R}^M$ we denote by $|a|_1 = \sum_{j=1}^{M} |a_j|$ the $\ell_1$ norm of a vector.

## 2. Sparse balls for the $\ell_1$ and $\ell_1 + \ell_2$ penalized estimates

In this section we establish upper bounds on the $\ell_1$ balls $|\widehat{\beta} - \beta^*|_1$ and $|\widetilde{\beta} - \beta^*|_1$, for the Lasso and elastic net estimates, in linear and logistic regression, respectively. We show that these bounds are, up to constants that we make precise below, of the form $k^* r$, where $r$ is the tuning parameter corresponding to the $\ell_1$ penalty and $k^*$ is the number of non-zero components of $\beta^*$. Since the $\ell_1$ norm is a sum of $M$ terms, but the bound only involves the unknown and possibly much lower dimension $k^*$, we call the corresponding balls sparse.

### 2.1. Conditions on the design matrix

In [5] and [19] it was showed that the Lasso type-estimates belong to sparse $\ell_1$ balls centered at the true parameter, in linear models and generalized linear models, respectively. These results were established under variants of a condition on the design matrix typically referred to as the mutual coherence condition, introduced in [8]. We state below a mild version of this condition, which we will also use in Section 3 of this paper. Let

$$\rho_{kj} = \frac{1}{n}\sum_{i=1}^{n} X_{ki} X_{ji}, \quad 1 \leq j, k \leq M.$$

*Condition Identif:*  We assume that there exists a constant $0 < d \leq 1$ such that

$$\mathbb{P}\left(\max_{j \in I^*, k \neq j} |\rho_{kj}| \leq \frac{d}{k^*}\right) = 1.$$

This condition guarantees separation of the variables in the true set $I^*$ from one another and from the rest, where the degree of separation is measured in terms of the size of the correlation coefficients. We regard it here as an identifiability condition. It will be used as a sufficient condition for correct variable selection in Section 3 below. However, it is not needed for sparse oracle inequalities, as we detail below.

In Sections 2.2 and 2.3 below we show that *Condition Identif* can be relaxed if one is only interested in prediction or the global behavior of the estimates measured, as in these sections, by the $\ell_1$ distance to the truth. To formulate the weaker condition let $\alpha > 0, \epsilon \geq 0$ be given. Define the set

$$V_{\alpha,\epsilon} = \left\{ v \in \mathbb{R}^M : \sum_{j \notin I^*} |v_j| \leq \alpha \sum_{j \in I^*} |v_j| + \epsilon \right\}. \tag{2.1}$$

Let $\Sigma$ be the $M \times M$ matrix with entries $\rho_{kj}$.

*Condition Stabil.* Let $\alpha, \epsilon > 0$ be given. There exist $0 < b \leq 1$ such that

$$\mathbb{P}\left( v'\Sigma v \geq b \sum_{j \in I^*} v_j^2 - \epsilon \right) = 1, \quad \text{for any } v \in V_{\alpha,\epsilon}.$$

*Remark.* We denote generically one of the estimates of $\beta^*$ studied below by $\breve{\beta}$. We will motivate the definition of the set $V_{\alpha,\epsilon}$ by showing, in the course of the proofs of Theorems 2.2–2.7, that $\breve{\beta} - \beta^* \in V_{\alpha,\epsilon}$, with high probability, for specific parameters $\alpha$ and $\epsilon$. For instance, we will show that $\alpha$ is either 3, for the $\ell_1$ penalized estimates, or 4, for the $\ell_1 + \ell_2$ penalized estimates. The parameter $\epsilon$ will be either zero, for the least squares estimates, or exponentially small, for each $M$ and $n$, in the case of the logistic regression estimates. The term $\epsilon$ in the definition of $V_{\alpha,\epsilon}$ is needed for purely technical reasons, and does not affect the results or their interpretation. *Condition Stabil* corresponding to $\alpha = 3$ and $\epsilon = 0$ has been introduced, for an analysis similar to the one we conduct here, by [2], for a comparative study of the predictive performance of the Dantzig and Lasso estimators in linear models.

One possible intuitive interpretation of *Condition Stabil* is as follows. If $\epsilon = 0$, *Condition Stabil* is immediately implied by $\mathbb{P}(\Sigma - b\mathrm{D} \geq 0) = 1$, where $D$ is the $M \times M$ matrix containing the $k^* \times k^*$ identity matrix corresponding to indices in $I^*$, and with zero elements otherwise. This asserts that the correlation matrix remains semi-positive definite if we decrease the diagonal elements corresponding to the true variables slightly, and leave all other entries unchanged. Since this modification affects only $k^*$ of $M^2$ entries, it can be regarded as a stability requirement on the correlation structure. *Condition Stabil* is even milder than $\mathbb{P}(\Sigma - b\mathrm{D} \geq 0) = 1$, since it is only required to hold for $v \in V_{\alpha,\epsilon}$, for some given $\alpha$ and $\epsilon$.

The following lemma establishes the relationship between the two conditions, and shows that *Condition Stabil* is less restrictive. A brief argument establishing this link is also offered in [2], for $\alpha = 3$ and $\epsilon = 0$; we include a full proof here for the general case, for completeness.

**Lemma 2.1.** *Let $\alpha > 0$ and $\epsilon \geq 0$ be given. If Condition Identif holds for some $0 < d < 1/(1+2\alpha+\epsilon)$, then Condition Stabil holds for any $0 < b \leq 1 - d(1+2\alpha+\epsilon)$.*

*Proof.* Let $\Sigma^*$ be the $k^* \times k^*$ matrix with entries $\rho_{kj}$, $k, j \in I^*$. For any $v \in \mathbb{R}^M$ denote by $v_*$ the vector in $\mathbb{R}^{k^*}$ obtained from $v$ by retaining only the components corresponding to $I^*$. Then

$$
\begin{aligned}
v'\Sigma v &\geq v'_* \Sigma^* v_* - 2 \sum_{j \in I^*} \sum_{k \notin I^*} |\rho_{kj}||v_j||v_k| \\
&\geq v'_* \Sigma^* v_* - \frac{2d}{k^*} \sum_{j \in I^*} |v_j| \sum_{k \notin I^*} |v_k|, \quad \text{under } Condition \ Identif \\
&\geq v'_* \Sigma^* v_* - \frac{2\alpha d}{k^*} \left( \sum_{j \in I^*} |v_j| \right)^2 - \frac{2d\epsilon}{k^*} \sum_{j \in I^*} v_j^2, \quad \text{for } v \in V_{\alpha,\epsilon} \\
&\geq v'_* \Sigma^* v_* - 2\alpha d \sum_{j \in I^*} v_j^2 - \frac{2d\epsilon}{k^*} \sqrt{k^*} \left( \sum_{j \in I^*} v_j^2 \right)^{1/2}, \quad \text{by Cauchy -Schwarz} \\
&\geq v'_* \Sigma^* v_* - (2\alpha + \epsilon) d \sum_{j \in I^*} v_j^2 - d\epsilon, \quad \text{since } 2xy \leq x^2 + y^2 \\
&\geq (1 - d(1 + 2\alpha + \epsilon)) \sum_{j \in I^*} v_j^2 - \epsilon.
\end{aligned}
$$

The last inequality follows from *Condition Identif*, which also implies that $v'_* \Sigma^* v_* \geq (1-d) \sum_{j \in I^*} v_j^2$ and so *Condition Stabil* holds for any $b$ with $0 < b \leq 1 - d(1 + 2\alpha + \epsilon)$. $\square$

Thus, for instance, for the study of the Lasso estimates in linear models, we have $\alpha = 3$ and $\epsilon = 0$ and so if *Condition Identif* holds for some $d$, then *Condition Stabil* holds for $0 < b \leq 1 - 7d$, which imposes the restriction $0 < d < \frac{1}{7}$.

The results of Sections 2.2 and 2.3 below will be established directly under the less restrictive *Condition Stabil*, which requires the specification of a constant $b$. Notice that if $b$ is very small, the condition is almost a tautology, as $\Sigma \geq 0$ by construction. However, as it will become apparent from the results established below, a very small value of $b$ will increase the radius of the $\ell_1$ balls covering the estimator. This motivates the parallel study of the elastic net estimates. We show that they are less affected by potentially small values of $b$.

## 2.2. *Sparse $\ell_1$ balls for estimates in linear regression models*

Throughout all sections on linear regression in this paper we assume that the model generating that data is $E(Y|\mathbf{X} = x) = \sum_{j \in I^*} \beta_j^* x_j$, for $\mathbf{X} \in \mathbb{R}^M$ and $I^* \subseteq \{1, \ldots, M\}$. This is the most popular model for regression with unbounded

response $Y$. It is also becoming increasingly common in regression models with $Y \in \{0,1\}$, when the data supports it. Its usage in this context dates back to [1].

### 2.2.1. An $\ell_1$ penalized least squares estimator

We estimate $\beta^*$ by

$$\widehat{\beta} = \arg\min_{\beta} \frac{1}{n} \sum_{i=1}^{n} \{Y_i - \beta' \mathbf{X}_i\}^2 + 2r \sum_{j=1}^{M} |\beta_j|, \qquad (2.2)$$

where $r =: r_{n,M}(\delta)$ is a tuning sequence depending on $n, M$ and a user specified parameter $\delta$. In what follows we determine $r$ such that $\mathbb{P}(|\widehat{\beta} - \beta^*|_1 \geq Crk^*) \geq 1 - \delta$, and we make $C > 0$ precise.

In the following theorem we will use *Condition Stabil* corresponding to the set $V_{\alpha,\epsilon}$ defined in (2.1) for $\epsilon = 0$ and $\alpha = 3$. Let $\sigma^2 = \text{Var}(Y)$ and recall that $L$ denotes a common bound on $X_j$, $1 \leq j \leq M$.

**Theorem 2.2.** *Assume that Condition Stabil is satisfied for some $0 < b \leq 1$. If we choose*

$$r \geq 2\sqrt{\frac{2\ln\frac{2M}{\delta}}{n}}, \; \text{if} \; Y \in \{0,1\},$$

*or*

$$r \geq 4L\sigma\sqrt{\frac{\ln\frac{4M}{\delta}}{n}} \vee 8L\frac{\ln\frac{4M}{\delta}}{n}, \; \text{if} \; Y \in \mathbb{R},$$

*then*

$$\mathbb{P}\left(|\widehat{\beta} - \beta^*|_1 \leq \frac{4}{b}rk^*\right) \geq 1 - \delta,$$

*for $\widehat{\beta}$ given in (2.2).*

*Remark 1.* In practice one can replace $\sigma$ in the tuning sequence by an estimator, as discussed in detail in [4].

*Remark 2.* It is interesting to note that although the results above indicate that the radius of the $\ell_1$ ball is small if $k^*r \leq 1$, the proofs make no use of this restriction on $k^*$; in particular $k^* > \sqrt{n}$ is allowed. It is clear that in this case the bounds are large but, perhaps surprisingly, this does not affect the validity of variable selection, for some design matrices. We discuss this in detail in the next section.

Theorem 2.2 above shows that the bound on $|\widehat{\beta} - \beta^*|_1$ becomes large if *Condition Stabil* is satisfied only for very small values of $b$. One remedy is provided by a slightly modified estimator, which retains the sparsity properties of the Lasso estimates, but is less affected by small values of $b$. The modified estimate will be penalized least squares with a combined $\ell_1$ and $\ell_2$ penalty and we discuss it in the next subsection.

### 2.2.2. An $\ell_1 + \ell_2$ penalized least squares estimator

We estimate now $\beta^*$ by $\widehat{\beta}$, where

$$\widehat{\beta} = \arg\min_{\beta} \frac{1}{n} \sum_{i=1}^{n} \{Y_i - \beta' \mathbf{X}_i\}^2 + 2r \sum_{j=1}^{M} |\beta_j| + c \sum_{j=1}^{M} \beta_j^2. \qquad (2.3)$$

As before, the goal is to find $r =: r_{n,M}(\delta)$ and $c =: c_{n,M}(\delta)$ for which we can construct sparse balls for the estimates.

In the following theorem we will use *Condition Stabil* corresponding to the the set $V_{\alpha,\epsilon}$ defined in (2.1) for $\epsilon = 0$ and $\alpha = 4$.

**Theorem 2.3.** *Assume that Condition Stabil is satisfied for some $0 < b \leq 1$. If $\max_{j \in I^*} |\beta_j^*| \leq B$, for some $B > 0$, independent of $n$, and if*

$$r \geq 2\sqrt{\frac{2\ln\frac{2M}{\delta}}{n}}, \quad c = \frac{r}{2B}, \quad if \quad Y \in \{0,1\},$$

*or*

$$r \geq 4L\sigma\sqrt{\frac{\ln\frac{4M}{\delta}}{n}} \vee 8L\frac{\ln\frac{4M}{\delta}}{n}, \quad c = \frac{r}{2B}, \quad if \quad Y \in \mathbb{R},$$

*then*

$$\mathbb{P}\left(|\widehat{\beta} - \beta^*|_1 \leq \frac{4.25}{b+c} r k^*\right) \geq 1 - \delta,$$

*for $\widehat{\beta}$ given in (2.3).*

*Remark.* The result above shows that even if *Condition Stabil* holds with $b$ very close to 0, the bound on $|\widehat{\beta} - \beta^*|_1$ stays finite, for any given $M$ and $n$. Note that it may still be large, since $c$ is restricted to take relatively small values, dictated by the sizes of $r$ and $B$. However, we cannot choose a much larger value for $c$: in that case the $\ell_2$ penalty would become prevalent, and no estimates will be set to zero in finite samples.

### 2.3. Sparse $\ell_1$ balls for estimates in logistic regression models

We denote the logistic loss function by

$$l(\beta) =: l(\beta; x, y) =: -y\beta'x + \log(1 + \exp\beta'x),$$

and denote by $\mathbb{P}l(\beta) = \mathbb{E}l(\beta; Y, \mathbf{X})$ the associated risk. Define

$$\beta^* = \arg\min_{\beta \in \mathbb{R}^M} \mathbb{P}l(\beta).$$

Throughout all sections on logistic regression we will assume that

$$\mathbb{E}(Y|\mathbf{X} = x) = p(x) = \frac{\exp(\sum_{j \in I^*} \beta_j^* x_j)}{1 + \exp(\sum_{j \in I^*} \beta_j^* x_j)}. \qquad (2.4)$$

### 2.3.1. An $\ell_1$ penalized logistic regression estimator

We estimate $\beta^*$ by

$$\widetilde{\beta} = \arg\min_{\beta} \frac{1}{n} \sum_{i=1}^{n} \left\{ -Y_i \sum_{j=1}^{M} \beta_j X_{ij} + \log\left(1 + \exp\sum_{j=1}^{M} \beta_j X_{ij}\right)\right\} + 2r\sum_{j=1}^{M} |\beta_j|.$$
(2.5)

We will determine the tuning sequence $r = r_{n,M}(\delta)$, different than the one above, for which we can construct sparse balls for these estimators. We will analyze the estimates under the assumption that $p(x)$ is bounded away from zero and one for all x. This is implied by:

*Assumption A*: There exists $D > 0$ such that $|\beta^*|_1 \le D$,

in combination with the assumption that all $X$ variables are bounded by $L$.

In the following theorem we will use *Condition Stabil* corresponding to the set $V_{\alpha,\epsilon}$ defined for

$$\epsilon = \frac{\log 2}{2^{(M \vee n)+1}} \times \frac{1}{r},$$

with $r$ given below, and for $\alpha = 3$. Also, let $s$ be a constant depending on $L$ and $D$, which decreases with $D$.

**Theorem 2.4.** *Assume Condition Stabil is satisfied for some $0 < b \le 1$. If Assumption A holds and if*

$$r \ge (6 + 4\sqrt{2})L\sqrt{\frac{2\log 2(M \vee n)}{n}} + 2L\sqrt{\frac{2\log(\frac{1}{\delta})}{n}} + \frac{1}{4(M \vee n)},$$

*then*

$$\mathbb{P}\left(|\widetilde{\beta} - \beta^*|_1 \le \frac{4}{sb}rk^* + (1 + \frac{1}{r})\epsilon\right) \ge 1 - \delta,$$
(2.6)

*for $\widetilde{\beta}$ given by (2.5).*

*Remark.* Notice that the term $(1 + \frac{1}{r})\epsilon$ is roughly $\frac{n}{2^{M \vee n}}$ and therefore negligible. As noted above, the bound on $|\widetilde{\beta} - \beta^*|_1$ becomes large for very small values of $b$, which motivates the study of the $\ell_1 + \ell_2$ penalized estimators in the next section. Also, the constant $1/s$, the exact form of which is given in the course of the proof of this theorem, can be very large for large $D$; similar results, based on different arguments and slightly more restrictive assumptions on $\Sigma$ have also been obtained in [19]. However, if we content ourselves with asymptotic statements, we can obtain an improved bound on $|\widetilde{\beta} - \beta^*|_1$, with $1/s$ replaced by a quantity arbitrarily close to 1. These results will hold with probability converging to 1, and under more stringent requirements on the design. They are based on the following fact, which is of independent interest: it establishes the sup-norm consistency of $\widetilde{\beta}'x$ as an estimate of $\beta^{*\prime}x$.

**Proposition 2.5.** *Let $\delta =: \delta_n$ be any sequence converging to zero with $n$. If $\max_{j \in I^*} |\beta_j^*| \le B$, for some $B > 0$, independent of $n$, and $k^* r \to 0$, then for any $\eta > 0$ we have*

$$\mathbb{P}\left(\sup_x |\widetilde{\beta}'x - \beta^{*'}x| \ge \eta\right) \longrightarrow 0,$$

*as $n \to \infty$.*

*Remark.* It is interesting to note that this result holds independently of the assumptions on the design.

The next result is obtained under a condition similar to *Condition Stabil*, but required to hold for a weighted version of the matrix $\Sigma$. Let $g(z) = e^z/(1+e^z)$. Let $\Sigma_1$ be the $M \times M$ matrix with entries $\frac{1}{n} \sum_{i=1}^n g'(\beta^{*'} X_i) X_{ki} X_{ji}$, $1 \le j, k \le M$.

*Condition LStabil.* Let $\alpha > 0, \epsilon \ge 0$ be given. There exist $0 < b \le 1$ such that

$$\mathbb{P}\left(v'\Sigma_1 v \ge b \sum_{j \in I^*} v_j^2 - \epsilon\right) = 1, \quad \text{for any} \ \ v \in V_{\alpha,\epsilon}.$$

**Theorem 2.6.** *Assume Condition LStabil is satisfied for some $0 < b \le 1$. Let*

$$r \ge (6 + 4\sqrt{2})L\sqrt{\frac{2\log 2(M \vee n)}{n}} + 2L\sqrt{\frac{2\log(\frac{1}{\delta_n})}{n}} + \frac{1}{4(M \vee n)},$$

*for any sequence $\delta_n \to 0$. If $\max_{j \in I^*} |\beta_j^*| \le B$, for some $B > 0$ independent of $n$, and $k^* r \to 0$, then*

$$\mathbb{P}\left(|\widetilde{\beta} - \beta^*|_1 \le \frac{4}{wb}rk^* + (1 + \frac{1}{r})\epsilon\right) \longrightarrow 1, \tag{2.7}$$

*for $\widetilde{\beta}$ given by (2.5) and for a constant $w$ arbitrarily close to one.*

### 2.3.2. *An $\ell_1 + \ell_2$ penalized logistic regression estimator*

In this section we obtain similar results for estimators of $\beta^*$ given by

$$\widetilde{\beta} = \arg\min_\beta \frac{1}{n}\left\{-Y_i \sum_{j=1}^M \beta_j X_{ij} + \log\left(1 + \exp\sum_{j=1}^M \beta_j X_{ij}\right)\right\}$$

$$+ 2r \sum_{j=1}^M |\beta_j| + c \sum_{j=1}^M \beta_j^2.$$

In the following theorem we will use *Condition Stabil* corresponding to the set $V_{\alpha,\epsilon}$ defined in (2.1) for

$$\epsilon = \frac{\log 2}{2^{(M \vee n)+1}} \times \frac{1}{r},$$

for $r$ given below, and $\alpha = 4$. Let $s > 0$ be the constant given in Theorem 2.4 above.

**Theorem 2.7.** *Assume that Assumption A holds and that Condition Stabil is satisfied for some $0 < b < 1$. Let $B > 0$ such that $\max |\beta|_j^* \le B$ and take*

$$r \ge (6 + 4\sqrt{2})L\sqrt{\frac{2\log 2(M \vee n)}{n}} + \frac{1}{4(M \vee n)} + 2L\sqrt{\frac{2\log(\frac{1}{\delta})}{n}}, \qquad c = \frac{r}{2B}.$$

*Then*

$$\mathbb{P}\left(|\widetilde{\beta} - \beta^*|_1 \le \frac{4.25 r k^*}{sb + c} + (1 + \frac{1}{r})\epsilon\right) \ge 1 - \delta,$$

*for $\widetilde{\beta}$ given by (2.8).*

The comments and remarks of Section 2.2.1 apply here with no change: the $\ell_1 + \ell_2$ penalized estimates are more stable, in that the radius of the $\ell_1$ ball covering the estimate is less affected by small values of $b$ and $s$. However, care must be exercised in choosing too large a $c$, as in that case the sparsity properties will be lost. We can also derive, in a similar manner, versions of Proposition 2.5 and Theorem 2.6 for the $\ell_1 + \ell_2$ penalized estimate. Since the results are nearly identical we do not include them here.

## 3. Correct subset selection

The asymptotic properties of subset selection via the Lasso in linear models have been studied recently by a number of authors: [13] studied selection Gaussian graphical models, [21] investigated subset selection in linear regression on for what was termed incoherent design matrices, [3] studied approximating regression models under design matrices satisfying *Condition Identif* introduced in Section 2 above and previously discussed in [4], and [25] investigated a three stage procedure in linear models. A nice overview of the connections between incoherent design matrices and matrices satisfying conditions similar to our *Condition Identif* is given in [14]. An interesting asymptotic analysis, in which one studies the interplay between the sample size $n$, the sparsity level $k^*$ and the number of variables $M$ for average asymptotic consistency in linear regression models with Gaussian design is presented in [24]. There the coefficient set $I^*$ is assumed to have been selected uniformly at random from $\{1, \ldots, M\}$, and one studies asymptotically the average error probability, where one averages over all possible choices of $I^*$. We refer to the work of [6] for a non-asymptotic investigation of the accuracy of model selection via the Lasso in linear models, but under model assumptions different than ours: the coefficient set $I^*$ is again assumed to have been selected uniformly at random from $\{1, \ldots, M\}$, and conditionally on $I^*$ the signs of $\beta_j$, $j \in I^*$, are assumed to be equally likely to be 1 or -1. The properties of the Lasso-type estimates used for correct subset selection in logistic regression have not been investigated from the perspective considered here. After finishing this paper, we learned of the recent work of [17], which investigates the very related topic of asymptotic model selection consistency in binary graphs; we comment on connections with our work in Section 3.2.2.

The finite sample properties of variable selection via the elastic net have not been investigated in either of the models considered here. For a discussion of its usage in linear regression models with different target parameters than those considered here we refer to [23].

We study in this section the non-asymptotic merits of the Lasso and elastic net estimates when used for variable selection. We conclude the section with the asymptotic implications of these results.

All estimates of $\beta^*$ analyzed in Theorems 2.2–2.7 have zero coefficients. These theorems, however, do not necessarily guarantee that the corresponding set of the non-zero coefficients of these estimates is exactly equal with $I^*$, with high probability: we can either omit some of the true variables or include variables that do not belong to $I^*$ while still being able to control the radii of the $\ell_1$ balls. In this section we find estimates of $\beta^*$ that have the properties discussed in Section 2 and for which, in addition, we have $\mathbb{P}(I^* = \hat{I}) \geq 1 - \gamma$, for some given small $\gamma > 0$. Since $\mathbb{P}(I^* = \hat{I}) \geq 1 - \mathbb{P}(I^* \not\subseteq \hat{I}) - \mathbb{P}(\hat{I} \not\subseteq I^*)$, we find the subset $\hat{I}$ such that

$$\mathbb{P}(I^* \subseteq \hat{I}) \geq 1 - \gamma_1 \quad \text{and} \quad \mathbb{P}(I^* \supseteq \hat{I}) \geq 1 - \gamma_2,$$

with $\gamma_1 + \gamma_2 = \gamma$.

### 3.1. Correct inclusion of all true variables in the selected set

In this section we discuss conditions under which we can obtain results of the type

$$\mathbb{P}(I^* \subseteq \hat{I}) \geq 1 - \gamma_1,$$

for some given $\gamma_1 > 0$, for estimates having the properties discussed in Section 2 above. Lemma 3.1 below shows what governs the size of $\mathbb{P}(I^* \subseteq \hat{I})$. We discuss in detail to which extent we can use the results of Section 2 directly for this study. Recall that the cardinality of $I^*$ is $k^*$.

**Lemma 3.1.** *Let $\beta^*$ and $\breve{\beta}$ be a combination parameter/estimator from Section 2. Let $\hat{I}$ be the index set of the non-zero components of $\breve{\beta}$. Then*

$$\mathbb{P}\left(I^* \not\subseteq \hat{I}\right) \leq \mathbb{P}\left(|\breve{\beta} - \beta^*|_1 \geq \min_{l \in I^*} |\beta_l^*|\right). \tag{3.1}$$

*Proof.* The following display follows directly from the definitions of $\hat{I}$ and $I^*$.

$$
\begin{aligned}
\mathbb{P}(I^* \not\subseteq \hat{I}) &\leq \mathbb{P}(j \notin \hat{I} \text{ for some } j \in I^*) \\
&\leq \mathbb{P}(\breve{\beta}_j = 0 \text{ and } \beta_j^* \neq 0, \text{ for some } j \in I^*) \\
&\leq \mathbb{P}(|\breve{\beta}_j - \beta_j^*| = |\beta_j^*|, \text{ for some } j \in I^*) \\
&\leq \mathbb{P}(|\breve{\beta}_j - \beta_j^*| \geq \min_{l \in I^*} |\beta_l^*|, \text{ for some } j \in I^*) \\
&\leq \mathbb{P}\left(|\breve{\beta} - \beta^*|_1 \geq \min_{l \in I^*} |\beta_l^*|\right).
\end{aligned}
$$

$\square$

### *3.1.1. Detection of large signals*

The purpose of this subsection is to point out that the study of $\mathbb{P}(I^* \subseteq \hat{I})$ via a direct application of the sparse oracle bounds derived in the previous section may lead to suboptimal results. We argue this in what follows. Inequality (4.15) in the proof of Lemma 3.1 above makes it clear that the rate at which $\mathbb{P}(I^* \not\subseteq \hat{I})$ decays is governed by how small we can make the probability of estimating a non-zero component of $\beta^*$ by zero. However, if we further bound this probability and arrive at (3.1), we can use Theorems 2.2–2.7 of the previous section directly. We thus arrive at the following corollary.

**Corollary 3.2.** *Let $0 < \delta < 1$ be fixed. Assume Condition Stabil holds, for the parameters specified in the statements of Theorems 2.2–2.7, respectively.*

1. *$\ell_1$ penalized least squares in linear models:*

   *If $\min_{j \in I^*} |\beta_j^*| \geq \frac{4}{b} r k^*$ with $r$ given by Theorem 2.2, then $\mathbb{P}(I^* \subseteq \hat{I}) \geq 1 - \delta$.*

2. *$\ell_1 + \ell_2$ penalized least squares in linear models:*

   *If $\frac{4.25}{b+c} r k^* \leq \min_{j \in I^*} |\beta_j^*| \leq \max_{j \in I^*} |\beta_j^*| \leq B$, for some $B > 0$, and with $r$ and $c$ given by Theorem 2.3, then $\mathbb{P}(I^* \subseteq \hat{I}) \geq 1 - \delta$.*

3. *$\ell_1$ penalized logistic regression:*

   *If $\min_{j \in I^*} |\beta_j^*| \geq \frac{4}{sb} r k^* + (1 + \frac{1}{r})\epsilon$, with $s$, $r$ and $\epsilon$ given by Theorem 2.4 and if Assumption A holds, then $\mathbb{P}(I^* \subseteq \hat{I}) \geq 1 - \delta$.*

4. *$\ell_1 + \ell_2$ penalized logistic regression:*

   *If $\frac{4.25}{sb+c} r k^* + (1 + \frac{1}{r})\epsilon \leq \min_{j \in I^*} |\beta_j^*| \leq \max_{j \in I^*} |\beta_j^*| \leq B$, with $r$, $c$ and $\epsilon$ given by Theorem 2.7 and if Assumption A holds, then $\mathbb{P}(I^* \subseteq \hat{I}) \geq 1 - \delta$.*

*Remark.* The lower bounds on the minimum size of the true coefficients stated in Corollary 3.2 are all of the type

$$\min_{j \in I^*} |\beta_j^*| \geq Cr k^*, \qquad (3.2)$$

possibly up to the small additive term $\epsilon$ defined in the previous section.

For stable design matrices, when the constant $C$ is close to 1, and if the true model is supported on a space of dimension $k^*$, with very low $k^*$ satisfying $r k^* < 1$, then such lower bounds imply that we can detect moderate sized signals. Clearly, for large $k^*$, the lower bounds on the coefficient size are too conservative, especially since the constant $C$ may also be large. We discuss below when one can weaken this requirement.

### 3.1.2. Detection of weaker signals

Propositions 3.3 and 3.4 below show that the lower bounds on the signal strength can be significantly weakened under further conditions on the design matrix. The intuition is the following: if the signal is very weak and the true variables are highly correlated with one another and with the rest, one cannot hope to recover the true model with high probability. We will therefore work, for the remainder of this paper, under the assumption that the true model is identifiable, as quantified in *Condition Identif* stated in Section 2 above. Recall that this condition only requires that the true variables be separated from on another and from the rest, and it does not impose any restrictions on the variables placed outside the true set.

**Detection of weak signals via $\ell_1$ and $\ell_1 + \ell_2$ penalized least squares in linear models.**

We show below that if the identifiability condition is met, then we can recover coefficients with sizes above the noise level $n^{-1/2}$. The following result shows that, if the identification is to be performed at some given confidence level $\delta$, the size of the signal will also depend on $\delta$. Moreover, it will depend on $M$, via a logarithmic term: this is the price to pay for simultaneous identification of the true variables, among all $M$ possibilities. In what follows we will use the following tuning parameters, depending whether $Y \in \{0,1\}$ or $Y \in \mathbb{R}$. Let $0 < \delta < 1$ be fixed. Let $K$ be an upper bound on $k^*$. Since $k^*$ is unknown, one can always use the conservative bound $M$. However, if in practical situations $K$ is known, one can use it instead of the larger bound $M$. Consider

$$
\begin{aligned}
r &\geq 2\sqrt{\frac{2\ln(\frac{2KM}{\delta})}{n}}, \text{ if } Y \in \{0,1\} \quad\quad (3.3)\\
&\quad\quad\quad\quad\text{or}\\
r &\geq 4L\sigma\sqrt{\frac{\ln\frac{4KM}{\delta}}{n}} \vee 8L\frac{\ln\frac{4KM}{\delta}}{n}, \text{ if } Y \in \mathbb{R}.
\end{aligned}
$$

**Proposition 3.3.** *For $r$ given above we assume that*

$$
\min_{j \in I^*} |\beta_j^*| \geq 2r.
$$

*(1) If Condition Identif is satisfied for $d \leq \frac{1}{15}$ and $\widehat{I}$ corresponds to the $\ell_1$ penalized least squares estimate, then*

$$
\mathbb{P}(I^* \subseteq \hat{I}) \geq 1 - \delta - \frac{\delta}{M}.
$$

*(2) We assume, in addition, that $\max_{j \in I^*} |\beta_j^*| \leq B$ for some $B > 0$. We choose $c = \frac{r}{2B}$. If Condition Identif is satisfied for $d \leq \frac{1+c}{17.5}$ and $\widehat{I}$ corresponds to the*

*$\ell_1 + \ell_2$ penalized least squares estimate, then*

$$\mathbb{P}(I^* \subseteq \hat{I}) \geq 1 - \delta - \frac{\delta}{M}.$$

*Remark.* Notice that although Corolarry 3.2 is established under the weaker *Condition Stabil*, it only guarantees that $\mathbb{P}(I^* \subseteq \hat{I})$ for the collection $I^*$ for which $\min_{j \in I^*} |\beta_j^*| \geq \frac{4}{b} r k^*$. In contrast, if *Condition Identif* holds, we can detect variables corresponding to the set $I^*$ for which $\min_{j \in I^*} |\beta_j^*| \geq 2r$. This is a substantial relaxation of the lower bound on the signal strength, which no longer depends on either the possibly large $k^*$ or the possibly small $b$. Similar relaxations of the requirements on $\min_{j \in I^*} |\beta_j^*|$ have also been obtained by [24] and [6], but for models in which $I^*$ is assumed to be random, as discussed at the beginning of Section 3.

Proposition 3.3 above allows an immediate comparison between the selection properties of the Lasso and the elastic net. Their behavior is almost the same, the only difference is in the restriction on the constant $d$: slightly larger values of $d$ can be allowed for the elastic net estimate. This translates into saying that if the correlations between the true variables, and between the true variables and the rest are slightly larger than what is allowed for the Lasso, then the $\ell_1 + \ell_2$ penalized estimate may provide an alternative. However, as we noted in Section 2, although it would be tempting to increase the value of $c$, in order to allow for a larger degree of correlation, this would result in not setting any components of the estimate to zero.

**Detection of weak signals via $\ell_1$ and $\ell_1 + \ell_2$ penalized logistic regression.**

The identifiability condition needed for linear models needs to be adjusted to the nature of the logistic regression model, in a manner similar to that of replacing *Condition Stabil* by *Condition LStabil*. We impose below a new condition: a weighted correlation matrix should exhibit the same type of separation we required of the correlation matrix of the data. The weights depend on the link function. This perhaps comes with little surprise: the correlation matrix appears explicitly in the expression of the least squares estimates in linear models, and this is not typically the case for other models and estimates. We formalize this below. For a given $0 < \delta < 1$, $M$ and $n$, let

$$r \geq (6 + 4\sqrt{2})L\sqrt{\frac{2\log 2(M \vee n)}{n}} + 2L\sqrt{\frac{2\log(\frac{2M}{\delta})}{n}} + \frac{1}{4(M \vee n)}, \quad (3.4)$$

$$\epsilon = \frac{\log 2}{2^{(M \vee n)+1}} \times \frac{1}{r},$$

where we recall that $L$ is a common bound on the $X_j$'s. Let $d$ be as required by *Condition Identif*. Recall that for such $0 < d < 1$ there exists a $0 < b < 1$ for which *Condition Stabil* holds, as specified in Lemma 2.1. For this $b$ define

$$U = \left\{ a \in \mathbb{R}^n : \max_{1 \le i \le n} \left| a_i - \sum_{j=1}^{M} \beta_j^* X_{ij} \right| \le \frac{4Lrk^*}{sb} + L(1 + \frac{1}{r})\epsilon \right\},$$

for $s > 0$ given in Theorem 2.4. The definition of $U$ is justified by the properties of the estimates $\widetilde{\beta}$ discussed in Section 2, which have been proved under *Condition Stabil* and *Assumption A*. Let $g(z) = e^z/(1 + e^z)$.

*Condition Lidentif.* Let $d$ be the constant required by *Condition Identif*. We assume that

$$\sup_{a \in U} \mathbb{P}\left( \max_{j \in I^*, k \ne j} \left| \frac{1}{n} \sum_{i=1}^n g'(a_i) X_{ij} X_{ik} \right| \le \frac{d}{k^*} \right) = 1.$$

*Remark 1.* We give a formal justification of this condition in the course of the proof of Proposition 3.4 in Appendix A below. It is a natural condition that appears via a linearization of the likelihood function. The term containing $\epsilon$ in the definition of $U$ is exponentially small, and can be essentially ignored for practical purposes; its role is purely technical.

**Proposition 3.4.** *Let $r$ and $\epsilon$ be as in (3.4) above and $s > 0$ as in Theorem 2.4. Let Assumption A hold.*

*(1) Assume that Conditions Identif and Lidentif are met with $d \le \frac{s}{16 + 2s(7+\epsilon)}$, for a set $U$ corresponding to $b \le 1 - d(7 + \epsilon)$. If*

$$\min_{j \in I^*} |\beta_j^*| \ge 3.5r + 3(1 + \frac{1}{r})\epsilon,$$

*and $\widehat{I}$ corresponds to the $\ell_1$ penalized logistic regression estimate then*

$$\mathbb{P}(I^* \subseteq \widehat{I}) \ge 1 - 3\delta.$$

*(2) Let $B > 0$ be such that $\max_{k \in I^*} |\beta_k^*| \le B$ and choose $c = \frac{2r}{B}$. Assume that Conditions Identif and Lidentif are met with $d \le \frac{s+c}{17 + 2s(8+\epsilon)}$, for a set $U$ corresponding to $b \le 1 - d(8 + \epsilon)$. If*

$$\min_{k \in I^*} |\beta_k^*| \ge 3.5r + (1 + \frac{1}{r})\epsilon,$$

*and $\widehat{I}$ corresponds to the $\ell_1 + \ell_2$ penalized logistic regression estimate then*

$$\mathbb{P}(I^* \subseteq \widehat{I}) \ge 1 - 3\delta.$$

*Remark 2.* Notice that if $g(x) = x$ is the linear link, *Condition Lidentif* becomes *Condition Identif*. Since $\epsilon$ is exponentially small, the requirement on the minimum size of the coefficients is essentially

$$\min_{k \in I^*} |\beta_k^*| \ge 3.5r.$$

As discussed in the remark following Proposition 3.3 above, Corollary 3.2 shows that $\mathbb{P}(I^* \subseteq \hat{I})$ can also be controlled under the less restrictive *Condition Stabil*, but in that case we can only recover sets $I^*$ corresponding to the large signal strength $\min_{j \in I^*} |\beta_j^*| \geq \frac{4}{sb} r k^* + (1 + \frac{1}{r}) \epsilon$. In contrast, Proposition 3.4 shows that we can detect weaker signals, however the correlation structure needs to follow the more restrictive *Conditions Lidentif* and *Identif*. As discussed before, similar properties are valid for the elastic net estimate, for an appropriate choice of the tuning sequence $c$. Refinements of this result, that replace the possibly small constant $s$ by a term close to 1 are possible, if instead of statements that hold with probability larger than $1 - \delta$ we consider statements that hold with probability converging to one. For this, one can use Proposition 2.5 and Theorem 2.6. Since these results are very similar to those above, we do not include them here, for brevity.

### 3.2. Correct subset selection

The set estimates $\hat{I}$ of the previous section satisfy $\mathbb{P}(I^* \subseteq \hat{I}) \geq 1 - \gamma_1$, for an appropriate $\gamma_1$. In what follows we show that $\hat{I}$ also satisfies $\mathbb{P}(\hat{I} \subseteq I^*) \geq 1 - \gamma_2$, thereby guaranteeing that $\mathbb{P}(\hat{I} = I^*) \geq 1 - \gamma$, for $\gamma_1 + \gamma_2 = \gamma$.

#### 3.2.1. Correct selection via the Lasso and the elastic net in linear regression models

**Theorem 3.5.** *Let $K$ be an upper bound on $k^*$ and take*

$$r \geq 2 \sqrt{\frac{2 \ln(\frac{2KM}{\delta})}{n}}, \; if \; Y \in \{0, 1\},$$

*or*

$$r \geq 4L\sigma \sqrt{\frac{\ln \frac{4KM}{\delta}}{n}} \vee 8L \frac{\ln \frac{4KM}{\delta}}{n}, \; if \; Y \in \mathbb{R}.$$

*Assume that*

$$\min_{j \in I^*} |\beta_j^*| \geq 2r$$

*(1) Assume that Condition Identif is met for $d \leq \frac{1}{15}$. If $\hat{I}$ corresponds to the $\ell_1$ penalized least squares estimator, then*

$$\mathbb{P}(\hat{I} = I^*) \geq 1 - 3\delta - \frac{\delta}{M}.$$

*(2) Assume, in addition, that $\max_{j \in I^*} |\beta_j^*| \leq B$ for some $B > 0$ and choose $c = \frac{r}{2B}$. If Condition Identif is met for $d = \frac{1+c}{17.5}$ and $\hat{I}$ corresponds to the $\ell_1 + \ell_2$ penalized least squares estimator, then*

$$\mathbb{P}(\hat{I} = I^*) \geq 1 - 3\delta - \frac{\delta}{M}.$$

*Remark.* Since $k^*$ is unknown, one can always take $K = M$. However, if in some instances one has a rough idea of the order of magnitude of $k^*$, one can use that value instead of the conservative bound $M$. The remarks on the relative merits of the Lasso versus the elastic net from the previous sections apply here with no change.

Recall that the Lasso parameter estimates $\widehat{\beta}$ may not be unique. However the set estimates $\hat{I}$ are unique, for each given tuning sequence $r$. This result, which we prove in Appendix B, is needed throughout the paper to ensure that the problem is well posed. We mention it again here, since it will be used constructively in the proof of Theorem 3.5 in Appendix A.

Theorem 3.5 has immediate asymptotic implications. It guarantees that $I^*$ will be consistently estimated by $\hat{I}$ if $M$, the number of candidate variables is polynomial in $n$, i.e $M = O(n^\zeta)$, for some $\zeta \geq 0$. To obtain this result it suffices to replace $\delta$ by any sequence converging to zero with $n$. For instance, choosing $\delta = 1/n$ and restating the value of $r$ in terms of order of magnitude we have the following corollary.

**Corollary 3.6.** *Let* $r = O(\sqrt{\frac{\log n}{n}})$ *and assume that* $\min_{j \in I^*} |\beta_j^*| = O(\sqrt{\frac{\log n}{n}})$. *Then, under the assumptions (1) or (2), respectively, of Theorem 3.5 we have*

$$\lim_{n \to \infty} \mathbb{P}(\hat{I} = I^*) = 1,$$

*for* $\hat{I}$ *either the* $\ell_1$ *or the* $\ell_1 + \ell_2$ *penalized least squares estimator.*

### 3.2.2. Correct variable selection via $\ell_1$ or $\ell_1 + \ell_2$ penalized logistic regression

In this subsection we show that the type of results that hold for $\ell_1$ or $\ell_1 + \ell_2$ penalized least squares continue to hold for penalized logistic regression, under requirements on the correlation matrix that are tailored to this type of loss function.

**Theorem 3.7.** *Under the assumptions of Proposition 3.4 we have:*

*(1) If* $\hat{I}$ *corresponds to the* $\ell_1$ *penalized logistic regression estimate then*

$$\mathbb{P}(I^* = \hat{I}) \geq 1 - 5\delta.$$

*(2) If* $\hat{I}$ *corresponds to the* $\ell_1 + \ell_2$ *penalized logistic regression estimate then*

$$\mathbb{P}(I^* = \hat{I}) \geq 1 - 5\delta.$$

The asymptotic implications of Theorem 3.7 are again immediate. If $M$ is polynomial in $n$ and for $\delta = 1/n$ we therefore obtain:

**Corollary 3.8.** *Let* $r = O(\sqrt{\frac{\log n}{n}})$ *and assume that* $\min_{j \in I^*} |\beta_j^*| = O(\sqrt{\frac{\log n}{n}})$. *Then, under Assumption A and the assumptions on the design required for (1) or (2), respectively, of Theorem 3.7 we have*

$$\lim_{n \to \infty} \mathbb{P}(\hat{I} = I^*) = 1,$$

*for $\widehat{I}$ either the $\ell_1$ or the $\ell_1 + \ell_2$ penalized logistic regression estimator.*

*Remark.* In the proofs of Corollary 3.8 and Theorem 3.7 above we invoked Theorem 2.4 of Section 2.3.1 and therefore used its hypotheses. If instead one invoked the asymptotic result of Theorem 2.6, one could obtain a version of Corollary 3.8 with Assumption A replaced by the conditions $\max_{k \in I^*} |\beta_k^*| \leq B$ and $rk^* \to 0$. For polynomially large $M$, the order of our tuning sequence is $r = O(\frac{\log n}{\sqrt{n}})$. The condition $rk^* \to 0$ therefore places the restriction $k^* \leq C \frac{\sqrt{n}}{(\log n)^2}$ on the size of the true model, for some positive constant $C$. In this context, similar results, based on different arguments, have been independently obtained by [17], under the slightly more stringent requirements $k^* \leq C(\frac{n}{\log n})^{1/3}$ and $\min_{k \in I^*} |\beta_k^*| \geq \frac{1}{k^*}$, but under slightly more relaxed conditions on the weighted matrix of the design.

## 4. Conclusions

The scope of this paper is to offer finite sample, non-asymptotic, benchmarks on the performance of the Lasso and the closely related elastic net methods for variable selection in logistic and linear regression methods. We showed that the methods can be used for correct variable selection in identifiable models, where we defined identifiability via *Condition Identif* and *Condition Lidentif*. The added requirement for correct selection, versus good prediction, is on the size of the signal strength: we can detect coefficients larger than a small constant multiplied by the tuning parameter of the $\ell_1$ penalty. This tuning parameter is a function of $n$, $M$ and the level of confidence, $\delta$. The size of the tuning parameter has to be larger than the noise level, typically of order $\frac{1}{\sqrt{n}}$, up to factors that are logarithmic in $M$ and $\frac{1}{\delta}$. Our contribution can be detailed as follows.

*Lasso and the elastic net in linear regression.* The properties of the $\ell_1$ penalized least squares in regression models are becoming well understood, while those of the $\ell_1 + \ell_2$ penalized least squares have not been investigated from this perspective. We complemented the existing results on the Lasso estimates by providing a refinement of assumptions. We showed in Section 2 that the $\ell_1$ penalized estimates belong to sparse $\ell_1$ balls under *Condition Stabil*, also proposed in [2]. We included a full proof of this result to facilitate the comparison with the elastic net estimates, which allow for a slightly higher degree of correlation between the **X** variables than the one permitted by the Lasso estimate. We discussed in Section 2 the precise interplay between this degree of correlation and the choice of the tuning parameters. If the tuning parameter of the $\ell_2$ term is smaller than the tuning parameter of the $\ell_1$ term, this estimator is also sparse: it belongs to a sparse $\ell_1$ ball centered at the true value and can be used to recover the true coefficient set $I^*$ with high probability. However, care must be taken when using this estimate: if the tuning sequence accompanying the extra $\ell_2$ term is too large we would essentially have a ridge regression estimate, and no variable selection will be performed.

In Section 3 we provided a non-asymptotic analysis of the subset selection problem in linear models, which complements the existing asymptotic results. We showed that the signal detection boundaries suggested by previous asymptotic analyses can be relaxed. In the works of [21] and [3], which investigate aspects of selection consistency, the minimal signal strength is required to be $n^{-\frac{1}{2}+\theta}$, for some $\theta > 0$, up to unspecified and possibly large constants. The work in [3] requires *Condition Identif* from Section 2 above. In [21] a less restrictive assumption on the design matrix is imposed, namely the irrepresentable design condition, which is almost necessary and sufficient for the sign consistency of the estimators, which implies consistent subset selection. The work of [14] uses a coherence-type condition similar to our *Condition Identif*, which is shown to be a sufficient condition for the sign consistency of a further thresholded Lasso estimator. The price to pay is a stronger requirement on the minimum size of the detectable coefficients: this size depends on sequences involved in the definition of their coherence condition and $k^*$. These requirements are similar in spirit to those discussed in our Corollary 3.2 above, and share similar drawbacks.

We showed here that if one concentrates directly on the study of $\mathbb{P}(\hat{I} = I^*)$, instead of sign consistency, and studies the original (untruncated) Lasso estimator under *Condition Identif*, one can relax the requirement on $\min_{j \in I^*} |\beta_j^*|$. We showed in Theorem 3.5 that one only needs $\min_{j \in I^*} |\beta_j^*|$ be larger than $\sqrt{\frac{2\ln(\frac{2M^2}{\delta})}{n}}$, up to small constants independent of the design. For $M$ polynomial in $n$ and the choice $\delta = 1/n$ one can therefore detect, with the untruncated Lasso, coefficients of order $O(\sqrt{\frac{\log n}{n}})$.

*Lasso and the elastic net in logistic regression models.* We showed in this article that the $\ell_1$ and $\ell_1 + \ell_2$ penalized logistic regression estimators have features that are similar to $\ell_1$ and $\ell_1 + \ell_2$ penalized least squares estimators, but the study of the estimates depends on conditions on a weighted correlation matrix of the data.

The predictive performance and adaptation to unknown sparsity of the Lasso penalized estimates in generalized linear models received very little attention, with the notable exceptions of [19, 26] and [11] in regression and classification. Here we revisited some of these issues, and showed that the $\ell_1$ penalized logistic regression estimators, as well as the elastic net estimates belong to sparse $\ell_1$ balls under *Condition Stabil*. The size of the radii of these balls can be improved asymptotically under *Condition LStabil*. We also showed that the $\ell_1 + \ell_2$ penalized logistic regression estimators, which have not yet been investigated, exhibit the same adaptation to unknown sparsity as the Lasso estimates, for appropriate choices of the tuning parameters given in Section 2.3. We showed in Theorem 3.7 that, similar to linear models, $\ell_1$ or $\ell_1 + \ell_2$ penalized logistic regression can be used to estimate $I^*$ with very high probability. The difference is in the conditions on the correlation matrix, which need to be adapted to the nature of this model, as in *Condition Lidentif*. The size of the coefficients that are detectable via this method is also of the order $O(\sqrt{\frac{\log n}{n}})$, where the constants involved in this bound are independent of the design or sparsity level.

## Appendix A

*Proof of Theorem 2.2.* Let $\mathbf{X}_i$ be the $M$ dimensional vector with entries $X_{ij}$, $1 \leq j \leq M$. For ease of notation, let $r_{n,M}(\delta) = r$. By the definition of the estimator, and with $W_i = Y_i - E(Y_i|\mathbf{X}_i)$, we obtain

$$(\beta^* - \widetilde{\beta})' \left\{ \frac{1}{n} \sum_{i=1}^{n} \mathbf{X}'_i \mathbf{X}_i \right\} (\beta^* - \widehat{\beta}) \tag{4.1}$$
$$\leq \sum_{j=1}^{M} |\beta_j^* - \widehat{\beta}_j| \left\{ \frac{2}{n} \left| \sum_{i=1}^{n} W_i X_{ij} \right| \right\} + 2r \sum_{j=1}^{M} |\beta_j^*| - 2r \sum_{j=1}^{M} |\widehat{\beta}_j|.$$

Define the event

$$A = \bigcap_{j=1}^{M} \left\{ \left| \frac{2}{n} \sum_{i=1}^{n} W_i X_{ij} \right| \leq r \right\}. \tag{4.2}$$

Notice that on the event $A$ display (4.1) yields, via simple algebra, that

$$\sum_{j \notin I^*} |\widehat{\beta}_j - \beta_j^*| \leq 3 \sum_{j \in I^*} |\widehat{\beta}_j - \beta_j^*|.$$

Therefore, on the set $A$ we have $\widehat{\beta} - \beta^* \in V$, with $V$ defined in (2.1), for $\epsilon = 0$ and $\alpha = 3$.

Adding $r|\widehat{\beta} - \beta^*|_1$ to both sides of (4.1) and re-arranging the terms we also have

$$(\beta^* - \widehat{\beta})' \left\{ \frac{1}{n} \sum_{i=1}^{n} \mathbf{X}'_i \mathbf{X}_i \right\} (\beta^* - \widehat{\beta}) + r|\widehat{\beta} - \beta^*|_1 \leq 4r \sum_{j \in I^*} |\beta_j^* - \widehat{\beta}_j|. \tag{4.3}$$

Using the Cauchy-Schwarz inequality in the right hand side of the inequality above, followed by an inequality of the type $2uv \leq au^2 + v^2/a$, for any $a > 1$, we further obtain

$$(\beta^* - \widehat{\beta})' \left\{ \frac{1}{n} \sum_{i=1}^{n} \mathbf{X}'_i \mathbf{X}_i \right\} (\beta^* - \widehat{\beta}) + r|\widehat{\beta} - \beta^*|_1 \leq 4ar^2 k^* + \frac{1}{a} \sum_{j \in I^*} (\beta_j^* - \widehat{\beta}_j)^2.$$

Since $\widehat{\beta} - \beta^* \in V$ we can invoke *Condition Stabil* and, by taking $a = 1/b$, we obtain, on the set $A$, that

$$|\widehat{\beta} - \beta^*|_1 \leq \frac{4}{b} r k^*.$$

To conclude the proof we determine now $r = r_{n,M}(\delta)$ such that $\mathbb{P}(A^c) \leq \delta$. If $Y \in \{0,1\}$ we use Hoeffding's inequality to obtain

$$\mathbb{P}(A^c) \leq \sum_{j=1}^{M} \mathbb{P}\left(\left|\frac{2}{n}\sum_{i=1}^{n} W_i X_{ij}\right| \geq r\right) \leq 2M \exp(-nr^2/8), \qquad (4.4)$$

and the choice

$$r_{n,M}(\delta) \geq 2\sqrt{\frac{2\ln(\frac{2M}{\delta})}{n}}$$

guarantees that $\mathbb{P}(A^c) \leq \delta$. If $Y \in \mathbb{R}$ we use Bernstein's inequality to obtain

$$
\begin{aligned}
\mathbb{P}(A^c) &\leq \sum_{j=1}^{M} \mathbb{P}\left(\left|\frac{2}{n}\sum_{i=1}^{n} W_i X_{ij}\right| \geq r\right) && (4.5) \\
&\leq 2M\left(\exp\left(-\frac{nr^2}{16L^2\sigma^2}\right) + \exp\left(-\frac{nr}{8L}\right)\right),
\end{aligned}
$$

and the choice

$$r_{n,M}(\delta) \geq 4L\sigma\sqrt{\frac{\ln\frac{4M}{\delta}}{n}} \vee \frac{8L}{n}\ln\frac{4M}{\sigma}$$

guarantees that $\mathbb{P}(A^c) \leq \delta$. This concludes the proof. $\qquad\square$

*Proof of Theorem 2.3.* Using the definition of the estimator, the fact that $\max_{j \in I^*}|\beta_j^*| \leq B$ and our choice of $c = \frac{r}{2B}$ we obtain, on the event $A$, that

$$\sum_{j \notin I^*}|\widehat{\beta}_j - \beta_j^*| \leq 4\sum_{j \in I^*}|\widehat{\beta}_j - \beta_j^*|.$$

Therefore, on the set $A$ we have $\widehat{\beta} - \beta^* \in V$, with $V$ defined in (2.1), for $\epsilon = 0$ and $\alpha = 4$.

We use the same reasoning as in Theorem 2.2, and invoke *Condition Stabil* to obtain the analogue of display (4.3). The only difference is that we complete the square generated by the $\ell_2$ part of the penalty:

$$
\begin{aligned}
b\sum_{j \in I^*}(\beta_j^* - \widehat{\beta}_j)^2 + c\sum_{j \in I^*}&(\beta_j^* - \widehat{\beta}_j)^2 + r|\widehat{\beta} - \beta^*|_1 && (4.6) \\
&\leq 2c\sum_{j \in I^*}\beta_j^*(\beta_j^* - \widehat{\beta}_j) + 4r\sum_{j \in I^*}|\beta_j^* - \widehat{\beta}_j|,
\end{aligned}
$$

and so, under the assumption that $\max_{j \in I^*}|\beta_j^*| \leq B$ and our choice of $c = \frac{r}{2B}$ we obtain, for any $a > 1$

$$(b+c)\sum_{j \in I^*}(\beta_j^* - \widehat{\beta}_j)^2 + r|\widehat{\beta} - \beta^*|_1 \leq 4.25ak^*r^2 + \frac{1}{a}\sum_{j \in I^*}(\beta_j^* - \widehat{\beta}_j)^2, \qquad (4.7)$$

and the remaining part of the proof is identical to that of Theorem 2.2, if we now choose $a = \frac{1}{b+c}$. $\qquad\square$

*Proof of Theorem 2.4.* Recall that we denoted the logistic loss function by

$$l(\beta) =: l(\beta; x, y) =: -y\beta'x + \log(1 + \exp\beta'x),$$

and the associated risk by $\mathbb{P}l(\beta) = \mathbb{E}l(\beta; Y, \mathbf{X})$. We also denote the empirical risk by

$$\mathbb{P}_n l(\beta) = \frac{1}{n} \sum_{i=1}^n \left\{ -Y_i \sum_{j=1}^M \beta_j X_{ij} + \log\left(1 + \exp \sum_{j=1}^M \beta_j X_{ij}\right) \right\}.$$

With this notation and letting $r = r_{n,M}(\delta)$, the estimator satisfies, by definition

$$\mathbb{P}_n l(\widetilde{\beta}) + 2r \sum_{j=1}^M |\widetilde{\beta}_j| \leq \mathbb{P}_n l(\beta^*) + 2r \sum_{j=1}^M |\beta_j^*|.$$

By adding and subtracting $\mathbb{P}(l(\widetilde{\beta}) - l(\beta^*)) + r \sum_{j=1}^M |\widetilde{\beta}_j - \beta_j^*|$ to both sides and rearranging terms we obtain

$$\begin{aligned}
r|\widetilde{\beta} - \beta^*|_1 + \mathbb{P}\left(l(\widetilde{\beta}) - l(\beta^*)\right) & \leq & (\mathbb{P}_n - \mathbb{P})\left(l(\beta^*) - l(\widetilde{\beta})\right) + r|\widetilde{\beta} - \beta^*|_1 \quad (4.8) \\
& & + \ 2r \sum_{j=1}^M |\beta_j^*| - 2r \sum_{j=1}^M |\widetilde{\beta}_j|.
\end{aligned}$$

Let

$$L_n = \sup_{\beta \in \mathbb{R}^M} \frac{(\mathbb{P}_n - \mathbb{P})(l(\beta^*) - l(\beta))}{|\beta - \beta^*|_1 + \epsilon}.$$

Notice first that if we change the $i$th pair $(\mathbf{X}_i, Y_i)$ while keeping the others fixed, the value of $L_n$ changes by at most $\frac{4L}{n}$, where $L$ is a common bound on all $X_{ij}$. To see why, recall that $\mathbb{P}_n = \frac{1}{n} \sum_{i=1}^n \delta_{\mathbf{X}_i, Y_i}$ is the empirical measure putting mass $1/n$ at each observation $(\mathbf{X}_i, Y_i)$. Let $\mathbb{P}'_n = \frac{1}{n}\left(\sum_{i=1, i \neq l}^{n-1} \delta_{\mathbf{X}_i, Y_i} + \delta_{\mathbf{X}'_l, Y'_l}\right)$ be the empirical measure corresponding to changing the pair $(\mathbf{X}_l, Y_l)$ to $(\mathbf{X}'_l, Y'_l)$. Then

$$\begin{aligned}
& \frac{(\mathbb{P}_n - \mathbb{P})(l(\beta^*) - l(\beta))}{|\beta - \beta^*|_1 + \epsilon} - \frac{(\mathbb{P}'_n - \mathbb{P})(l(\beta^*) - l(\beta))}{|\beta - \beta^*|_1 + \epsilon} \\
& = \frac{1}{n} \frac{l(\beta^*; Y_l, \mathbf{X}_l) - l(\beta; Y_l, \mathbf{X}_l) - l(\beta^*; Y'_l, \mathbf{X}'_l) + l(\beta; Y'_l, \mathbf{X}'_l)}{|\beta - \beta^*|_1 + \epsilon} \\
& \leq \frac{4L}{n} \frac{|\beta - \beta^*|_1}{|\beta - \beta^*|_1 + \epsilon} \leq \frac{4L}{n}, \qquad\qquad\qquad (4.9)
\end{aligned}$$

where the inequality follows immediately by a first order Taylor expansion and the assumption that all $X$ variables are bounded by $L$. Therefore we can apply the bounded difference inequality (e.g. Theorem 2.2, page 8 in [7]) to obtain that

$$\mathbb{P}(L_n - \mathbb{E}L_n \geq u) \leq \exp{-\frac{nu^2}{8L^2}}.$$

Thus, if we take

$$u > 2L\sqrt{\frac{2\log\frac{1}{\delta}}{n}},$$

we have

$$\mathbb{P}(L_n - \mathbb{E}L_n \geq u) \leq \delta.$$

We will use Lemma 3 in [26] to obtain a bound on $\mathbb{E}L_n$. We re-state it here for ease of reference, adapting it to our notation.

*Lemma 3, [26]. Let $J_n$ be an integer such that $2^{J_n} \geq n$. Then, if $L_n$ defined above corresponds to a Lipschitz loss and the components of $\mathbf{X}$ are bounded by $L$, with probability one, then $\mathbb{E}L_n \leq C_1\sqrt{\frac{2\log 2(M \vee n)}{n}} + C_2 \frac{J_n}{2^{(M \vee n)^2}}$, where $C_1, C_2$ are positive constants depending on the Lipschitz constant and $L$.*

Our loss is Lipschitz in $t = \beta'x$, with constant 2. Also, inspection of the chaining argument used in the proof of the Lemma shows that we can take $J_n = (M \vee n)$ and $\epsilon = \frac{\log 2}{2^{(M \vee n)+1}} \times \frac{1}{r}$. Therefore, by making the constants precise we obtain

$$\mathbb{E}L_n \leq 6L\sqrt{\frac{2\log 2(M \vee n)}{n}} + \frac{1}{4(M \vee n)}.$$

Define the event

$$E = \{L_n \leq r\}. \tag{4.10}$$

From the previous displays we then conclude that if

$$r \geq 6L\sqrt{\frac{2\log 2(M \vee n)}{n}} + \frac{1}{4(M \vee n)} + 2L\sqrt{\frac{2\log(\frac{1}{\delta})}{n}},$$

then $P(E) \geq 1 - \delta$.

Since $\mathbb{P}\left(l(\widetilde{\beta}) - l(\beta^*)\right) > 0$, by the definition of $\beta^*$, display (4.8) yields

$$
\begin{aligned}
r\sum_{j=1}^{M}|\widetilde{\beta}_j - \beta_j^*| &\leq \sum_{j=1}^{M} 2r|\widetilde{\beta}_j - \beta_j^*| + \sum_{j=1}^{M} 2r|\beta_j^*| - \sum_{j=1}^{M} 2r|\widetilde{\beta}_j| + r\epsilon \\
&\leq 2r\sum_{j=1}^{M}|\widetilde{\beta}_j| + 2r\sum_{j=1}^{M}|\beta_j^*| + 2r\sum_{j=1}^{M}|\beta_j^*| - \sum_{j=1}^{M} 2r|\widetilde{\beta}_j| + r\epsilon \\
&\leq 4r|\beta^*|_1 + r\epsilon,
\end{aligned}
$$

on the set $E$. Therefore, if Assumption A holds, we obtain

$$\sum_{j=1}^{M}|\widetilde{\beta}_j - \beta_j^*| \leq 4D + \epsilon \leq 5D, \tag{4.11}$$

where we used the possibly conservative bound $D$ on $\epsilon$ to keep the exposition clear.

By Example 4.5 in [18] we have $\mathbb{P}l(\widetilde{\beta}) - \mathbb{P}l(\beta^*) \geq \|g_{\widetilde{\beta}} - g_{\beta^*}\|^2$, where $g_\beta(x) = \frac{\exp(\beta'x)}{1+\exp(\beta'x)}$ and $\|\ \|$ is the $L_2$ norm with respect to the distribution of $\mathbf{X}$. A first order Taylor expansion gives $g_{\widetilde{\beta}}(x) - g_{\beta^*}(x) = \frac{\exp(\bar{\beta}'x)}{(1+\exp(\bar{\beta}'x))^2}(f_{\widetilde{\beta}}(x) - f_{\beta^*}(x))$, where $f_\beta(x) = \beta'x$ and $\bar{\beta}'x$ is an intermediate point between $\widetilde{\beta}'x$ and $\beta^{*'}x$. Let $A = 6LD$, and let $s = (1 + e^A)^{-4}$. Then, since Assumption A and (4.11) hold, we have $\|g_{\widetilde{\beta}} - g_{\beta^*}\|^2 \geq s\|f_{\widetilde{\beta}} - f_{\beta^*}\|^2$.

Thus, on the event $E$, display (4.8) further yields

$$
\begin{aligned}
r\sum_{j=1}^{M}|\widetilde{\beta}_j - \beta_j^*| + s\|f_{\widetilde{\beta}} - f_{\beta^*}\|^2 \quad &\leq \quad \sum_{j=1}^{M} 2r|\widetilde{\beta}_j - \beta_j^*| \qquad\qquad (4.12) \\
&+ \quad \sum_{j=1}^{M} 2r|\beta_j^*| - \sum_{j=1}^{M} 2r|\widetilde{\beta}_j| + r\epsilon \\
&\leq \quad 4\sum_{j\in I^*} r|\widetilde{\beta}_j - \beta_j^*| + r\epsilon.
\end{aligned}
$$

Via simple algebra, display (4.12) yields

$$
\sum_{j\notin I^*}|\widetilde{\beta}_j - \beta_j^*| \leq 3\sum_{j\in I^*}|\widetilde{\beta}_j - \beta_j^*| + \epsilon,
$$

on the set $E$. Therefore $\widetilde{\beta} - \beta^* \in V$, for the set $V$ given by (2.1) of Section 2, with $\alpha = 3$ and $\epsilon$ as in the statement of the theorem.

Let $\gamma_{kj} = \mathbb{E}X_k X_j$, for $k, j \in \{1, \ldots, M\}$ and let $\Gamma$ be the $M \times M$ matrix with entries $\gamma_{kj}$. Notice that $\|f_{\widetilde{\beta}} - f_{\beta^*}\|^2 = (\beta^* - \widetilde{\beta})'\Gamma(\beta^* - \widetilde{\beta})$ and so, using a reasoning identical to the one used in display (4.3) of Theorem 2.2, we further obtain

$$
r\sum_{j=1}^{M}|\widetilde{\beta}_j - \beta_j^*| + s(\beta^* - \widetilde{\beta})'\Gamma(\beta^* - \widetilde{\beta}) \leq 4ar^2k^* + \frac{1}{a}\sum_{j\in I^*}(\beta_j^* - \widetilde{\beta}_j)^2 + r\epsilon.
$$

Since on the set $E$ we have $\widetilde{\beta} - \beta^* \in V$, we can use *Condition Stabil*. The condition implies that $(\beta^* - \widetilde{\beta})'\Gamma(\beta^* - \widetilde{\beta}) \geq sb\sum_{j\in I^*}(\beta_j^* - \widetilde{\beta}_j)^2 - \epsilon$ and so

$$
r\sum_{j=1}^{M}|\widetilde{\beta}_j - \beta_j^*| + sb\sum_{j\in I^*}(\beta_j^* - \widetilde{\beta}_j)^2 \leq 4ar^2k^* + \frac{1}{a}\sum_{j\in I^*}(\beta_j^* - \widetilde{\beta}_j)^2 + (r+1)\epsilon. \quad (4.13)
$$

Taking $a = 1/sb$ we obtain, on the set $E$, that

$$
\sum_{j=1}^{M}|\widetilde{\beta}_j - \beta_j^*| \leq \frac{4rk^*}{sb} + (1 + \frac{1}{r})\epsilon.
$$

Since we have shown above that $P(E) \geq 1 - \delta$, the proof is complete. $\qquad\square$

*Proof of Proposition 2.5.* First notice that on the event $E$ defined in the previous theorem, display (4.8) yields

$$\mathbb{P}\left(l(\widetilde{\beta}) - l(\beta^*)\right) \quad \leq \quad 4r|\beta^*|_1 + r\epsilon.$$

Then, the assumptions of this proposition imply that the righthandside of the above display converges to zero with $n$ and so we have that for any $\vartheta > 0$

$$\mathbb{P}\left(|\mathbb{P}(l(\widetilde{\beta}) - l(\beta^*)| \geq \vartheta\right) \to 0, \tag{4.14}$$

since $\mathbb{P}(E^c) \leq \delta_n \to 0$.

Observe that $\mathbb{P}l(\beta) = \mathbb{P}_{\mathbf{X}}\mathbb{P}_{Y|\mathbf{X}}l(\beta)$, where we regard the expectations as being taken with respect to a pair $(\mathbf{X}, Y)$ independent of the sample. By the definitions of the loss $l$ and $p(x)$ we have

$$\mathbb{P}_{Y|\mathbf{X}=x}l(\beta) = \log(1 + e^{\beta'x}) - p(x)\beta'x.$$

Let $\theta > 0$ be arbitrary, fixed. Simple algebra shows that if $\sup_x |\widetilde{\beta}'x - \beta^{*'}x| \geq \theta$ then $\mathbb{P}_{Y|X=x}l(\widetilde{\beta}) > \mathbb{P}_{Y|X=x}l(\beta^*)$, for all $x$, and so there exists $\vartheta_\theta > 0$ such that $\mathbb{P}(l(\widetilde{\beta}) - l(\beta^*)) \geq \vartheta_\theta$. Then

$$\mathbb{P}\left(\sup_x |\widetilde{\beta}'x - \beta^{*'}x| \geq \theta\right) \leq \mathbb{P}\left(\mathbb{P}(l(\widetilde{\beta}) - l(\beta^*)) \geq \vartheta_\theta\right) \longrightarrow 0,$$

where the convergence to zero follows by (4.14) above. This concludes the proof of this proposition. □

*Proof of Theorem 2.6.* The proof differs from the proof of Theorem 2.4 above only in the way we obtain the lower bound on $\mathbb{P}(l(\widetilde{\beta}) - l(\beta^*))$. For quantities defined in the discussion immediately following display (4.11) above we write

$$\frac{\exp(\bar{\beta}'x)}{(1 + \exp(\bar{\beta}'x))^2}\left(f_{\widetilde{\beta}}(x) - f_{\beta^*}(x)\right)$$

$$= \exp(\bar{\beta}'x - \beta^{*'}x)\left(\frac{1 + \exp(\beta^{*'}x)}{1 + \exp(\bar{\beta}'x)}\right)^2(\widetilde{\beta} - \beta^*)'xp(x)(1 - p(x)),$$

where we recall that $\bar{\beta}'x$ is an intermediate point between $\widetilde{\beta}'x$ and $\beta^{*'}x$ and that we defined $p(x) = \frac{\exp\beta^{*'}x}{1 + \exp(\beta^{*'}x)}$. Let $\theta > 0$ be arbitrarily close to zero. Let $A_\theta$ be the set for which

$$\sup_x |\widetilde{\beta}'x - \beta^{*'}x| \leq \theta,$$

and recall that, by Proposition 2.5 we have $\mathbb{P}(A_\theta) \to 1$. On the set $A_\theta$ we have

$$\exp(\bar{\beta}'x - \beta^{*'}x)\left(\frac{1 + \exp(\beta^{*'}x)}{1 + \exp(\bar{\beta}'x)}\right)^2 \geq e^{-2\theta} =: w,$$

for all $x$ and with $w$ arbitrarily close to 1. Let $\Gamma_1$ be the matrix with entries $\mathbb{E}p(X)(1-p(X))X_kX_j$, for $k,j \in \{1, \ldots, M\}$. Therefore, on $A_\theta$ we have $\|g_{\widetilde{\beta}} - g_{\beta^*}\|^2 \geq w(\beta^* - \widetilde{\beta})'\Gamma_1(\beta^* - \widetilde{\beta})$. Invoking condition *LStabil* we obtain $(\beta^* - \widetilde{\beta})'\Gamma_1(\beta^* - \widetilde{\beta}) \geq wb\sum_{j \in I^*}(\widetilde{\beta}_j - \beta_j^*)^2$ and the rest of the proof carries on unchanged, with results holding now on the set $A_\theta \cap E$. $\qquad\square$

*Proof of Theorem 2.7.* The proof is identical to the one of Theorem 2.4 above, up to the following display

$$r\sum_{j=1}^{M}|\widetilde{\beta}_j - \beta_j^*| + sb\sum_{j \in I^*}(\beta_j^* - \widetilde{\beta}_j)^2 + c\sum_{j=1}^{M}\widetilde{\beta}_j^2 - c\sum_{j=1}^{M}\beta_j^{*2}$$
$$\leq 4r\sum_{j \in I^*}|\beta_j^* - \widetilde{\beta}_j| + (r+1)\epsilon.$$

To arrive at this display we observe that the elastic net satisfies

$$\sum_{j \notin I^*}|\widetilde{\beta}_j - \beta_j^*| \leq 4\sum_{j \in I^*}|\widetilde{\beta}_j - \beta_j^*| + \epsilon,$$

and so $\widetilde{\beta} - \beta^* \in V$, for the set $V$ given by (2.1) of Section 2, with $\alpha = 4$ and $\epsilon$ as in the statement of the theorem. Therefore the use of *Condition Stabil* in the derivations above is valid.

For the remaining of the proof we reason as in Theorem 2.3 above. We complete the square in the left hand side of the inequality above and invoke the assumption $\max_{j \in I^*}|\beta_j^*| \leq B$ to obtain

$$r\sum_{j=1}^{M}|\widetilde{\beta}_j - \beta_j^*| + sb\sum_{j \in I^*}(\beta_j^* - \widetilde{\beta}_j)^2 + c\sum_{j \in I^*}(\widetilde{\beta}_j - \beta_j^*)^2$$
$$\leq 4r\sum_{j \in I^*}|\beta_j^* - \widetilde{\beta}_j| + (r+1)\epsilon + 2cB\sum_{j \in I^*}|\beta_j^* - \widetilde{\beta}_j|,$$

which immediately implies, by choosing $c$ such that $2cB = r$, that

$$r\sum_{j=1}^{M}|\widetilde{\beta}_j - \beta_j^*| + (sb+c)\sum_{j \in I^*}(\beta_j^* - \widetilde{\beta}_j)^2 \quad \leq \quad 2 \times 2.5r\sum_{j \in I^*}|\beta_j^* - \widetilde{\beta}_j| + (r+1)\epsilon.$$

Then, we use again the Cauchy-Schwarz inequality followed by $2xy \leq ax^2 + y^2/a$ to obtain

$$r\sum_{j=1}^{M}|\widetilde{\beta}_j - \beta_j^*| + (sb+c)\sum_{j \in I^*}(\beta_j^* - \widetilde{\beta}_j)^2$$
$$\leq (2.5)^2ak^*r^2 + \frac{1}{a}\sum_{j \in I^*}(\beta_j^* - \widetilde{\beta}_j)^2 + (r+1)\epsilon.$$

Choosing now $a = \frac{1}{sb+c}$ gives, on the event $E$ defined in (4.10) above

$$\sum_{j=1}^{M} |\widetilde{\beta}_j - \beta_j^*| \leq \frac{4.25r}{sb+c} k^* + \left(1 + \frac{1}{r}\right)\epsilon.$$

Since we showed in the proof of Theorem 2.4 that $\mathbb{P}(E^c) \leq \delta$, this completes the proof. $\qquad\square$

*Proof of Proposition 3.3.* Recall that we denoted the cardinality of $I^*$ by $k^*$. First observe that by the definitions of $\widehat{I}$ and $I^*$ and by the union bound we have

$$
\begin{aligned}
\mathbb{P}(I^* \not\subseteq \widehat{I}) &\leq \mathbb{P}\left(k \notin \widehat{I} \text{ for some } k \in I^*\right) \\
&\leq \mathbb{P}\left(\widehat{\beta}_k = 0 \text{ and } \beta_k^* \neq 0, \text{ for some } k \in I^*\right) \\
&\leq k^* \max_{k \in I^*} \mathbb{P}\left(\widehat{\beta}_k = 0 \text{ and } \beta_k^* \neq 0\right).
\end{aligned}
$$

We first show that $\mathbb{P}(I^* \subseteq \widehat{I}) \geq 1 - \delta - \frac{\delta}{M}$ for the $\ell_1$ penalized least squares estimator. It follows immediately from Lemma 4.1 in Appendix B below that if $\widehat{\beta}_k = 0$ is a component of the solution $\widehat{\beta}$ then

$$\left|\frac{2}{n}\sum_{i=1}^{n}\left(Y_i - \sum_{j=1}^{M}\widehat{\beta}_j X_{ij}\right)X_{ik}\right| \leq 2r.$$

Therefore

$$
\begin{aligned}
&\mathbb{P}(I^* \not\subseteq \widehat{I}) \\
&\leq k^* \max_{k \in I^*} \mathbb{P}\left(\widehat{\beta}_k = 0 \text{ and } \beta_k^* \neq 0\right) \\
&\leq k^* \max_{k \in I^*} \mathbb{P}\left(\left|\frac{2}{n}\sum_{i=1}^{n}\left[Y_i - \sum_{j=1}^{M}\widehat{\beta}_j X_{ij}\right]X_{ik}\right| \leq 2r; \quad \beta_k^* \neq 0\right) \\
&= k^* \max_{k \in I^*} \mathbb{P}\left(\left|2\beta_k^* + \frac{2}{n}\sum_{i=1}^{n}W_i X_{ik} + \sum_{j \neq k}(\widehat{\beta}_j - \beta_j^*)\left(\frac{2}{n}\sum_{i=1}^{n}X_{ij}X_{ik}\right)\right| \leq 2r\right) \\
&\leq k^* \max_{k \in I^*} \mathbb{P}\left(|\beta_k^*| - \left|\frac{1}{n}\sum_{i=1}^{n}W_i X_{ik}\right| - \left|\sum_{j \neq k}(\widehat{\beta}_j - \beta_j^*)\left(\frac{1}{n}\sum_{i=1}^{n}X_{ij}X_{ik}\right)\right| \leq r\right) \\
&= k^* \max_{k \in I^*} \mathbb{P}\left(\left|\frac{1}{n}\sum_{i=1}^{n}W_i X_{ik}\right| + \sum_{j \neq k}|\widehat{\beta}_j - \beta_j^*|\left|\frac{1}{n}\sum_{i=1}^{n}X_{ij}X_{ik}\right| \geq |\beta_k^*| - r\right),
\end{aligned}
$$

where the penultimate inequality follows by the triangle inequality $|a+b+c| \geq |c| - |a| - |b|$. Under *Condition Identif* and since $\min_{j \in I^*} |\beta_j^*| \geq 2r$ we further

obtain

$$\mathbb{P}(I^* \nsubseteq \hat{I}) \quad \leq \quad k^* \max_{k \in I^*} \mathbb{P}\left( \left| \frac{1}{n} \sum_{i=1}^{n} W_i X_{ik} \right| \geq r/2 \right) + k^* \mathbb{P}\left( |\hat{\beta} - \beta^*|_1 \geq \frac{rk^*}{2d} \right).$$

We argue exactly as in the course of the proof of Theorem 2.2 to bound the probabilities above. We use either Hoeffding's inequality, for $Y \in \{0, 1\}$ or Bernstein's inequality, for $Y \in \mathbb{R}$ to bound the first term by $\frac{k^* \delta}{KM} \leq \frac{\delta}{M}$, for $r$ given by (3.3). Similarly, for this choice of $r$, we have

$$k^* \mathbb{P}\left( |\hat{\beta} - \beta^*|_1 \geq \frac{4rk^*}{b} \right) \leq k^* \times \frac{\delta}{K} \leq \delta,$$

for a constant $b$ for which *Condition Stabil* holds. By Lemma 2.1 in Section 2, *Condition Identif* implies *Condition Stabil* with $b = 1 - 7d$. Notice that for this value of $b$ we have $1/2d \geq 4/b$ for $d \leq 1/15$, as required in statement of this theorem. Therefore, combining these results we obtain

$$\mathbb{P}(I^* \nsubseteq \hat{I}) \quad \leq \quad \frac{\delta}{M} + \delta.$$

We establish now similar results for the $\ell_1 + \ell_2$ penalized least squares estimator. By the characterization of a zero component of the solution, given in Lemma 4.3 in Appendix B below, we also have

$$\mathbb{P}(I^* \nsubseteq \hat{I}) \leq \mathbb{P}\left( \left| \frac{2}{n} \sum_{i=1}^{n} [Y_i - \sum_{j=1}^{M} \hat{\beta}_j X_{ij}] X_{ik} \right| \leq 2r, \quad \beta_k^* \neq 0 \right),$$

and so the proof is identical to the one above. The only modification is in terms of constants: in this case *Condition Identif* implies *Condition Stabil* with $b = 1 - 9d$. From Theorem 2.3 we obtain for the choice of $r$ given by (3.3) that

$$\mathbb{P}\left( |\hat{\beta} - \beta^*|_1 \geq \frac{4.25rk^*}{b+c} \right) \leq \frac{\delta}{K}.$$

As above, we note $1/2d \geq 4.25/(b+c)$ for $d \leq \frac{1+c}{17.5}$. Invoking now Theorem 2.3 with these constants concludes the proof. □

*Proof of Proposition 3.4.* As in the previous proof, recall that we denoted the cardinality of $I^*$ by $k^*$ and that

$$\mathbb{P}(I^* \nsubseteq \hat{I}) \quad \leq \quad \mathbb{P}\left( k \notin \hat{I} \text{ for some } k \in I^* \right)$$
$$\leq \quad k^* \max_{k \in I^*} \mathbb{P}\left( \hat{\beta}_k = 0 \text{ and } \beta_k^* \neq 0 \right).$$

We begin by establishing the result for the $\ell_1$ penalized estimator. By Lemma 4.1 in the Appendix below it follows that if $\tilde{\beta}_k = 0$ is a component of the solution $\tilde{\beta}$ then

$$\left| \frac{1}{n} \sum_{i=1}^{n} X_{ik} \frac{\exp \sum_{j=1}^{M} \tilde{\beta}_j X_{ij}}{1 + \exp \sum_{j=1}^{M} \tilde{\beta}_j X_{ij}} - \frac{1}{n} \sum_{i=1}^{n} Y_i X_{ik} \right| \leq 2r.$$

Let now

$$S_n = \frac{1}{n} \sum_{i=1}^{n} X_{ik} \left\{ \frac{\exp \sum_{j=1}^{M} \widetilde{\beta}_j X_{ij}}{1 + \exp \sum_{j=1}^{M} \widetilde{\beta}_j X_{ij}} - \frac{\exp \sum_{j=1}^{M} \beta_j^* X_{ij}}{1 + \exp \sum_{j=1}^{M} \beta_j^* X_{ij}} \right\}.$$

Then, since $Y_i = Y_i - p(X_i) + p(X_i) =: W_i + p(X_i)$, where $p(X_i)$ is given by (2.4), we obtain:

$$\mathbb{P}(I^* \not\subseteq \hat{I}) \leq k^* \max_{k \in I^*} \mathbb{P}\left( \left| S_n - \frac{1}{n} \sum_{i=1}^{n} W_i X_{ik} \right| \leq 2r; \ \ \beta_k^* \neq 0 \right).$$

Define

$$B_n = \sum_{j=1}^{M} (\widetilde{\beta}_j - \beta_j^*) \left( \frac{1}{n} \sum_{i=1}^{n} X_{ij} X_{ik} \right).$$

Recalling that $\frac{1}{n} \sum_{i=1}^{n} X_{ik}^2 = 1$ we obtain, for every $k \in I^*$, that

$$\mathbb{P}\left( \left| S_n - B_n + B_n - \frac{1}{n} \sum_{i=1}^{n} W_i X_{ik} \right| \leq 2r; \ \ \beta_k^* \neq 0 \right)$$

$$\leq \mathbb{P}\left( |B_n| - |S_n - B_n| - \left| \frac{1}{n} \sum_{i=1}^{n} W_i X_{ik} \right| \leq 2r; \ \ \beta_k^* \neq 0 \right)$$

$$= \mathbb{P}\left( |\beta_k^*| - \left| \sum_{j \neq k} (\widetilde{\beta}_j - \beta_j^*) \frac{1}{n} \sum_{i=1}^{n} X_{ij} X_{ik} \right| - |S_n - B_n| - \left| \frac{1}{n} \sum_{i=1}^{n} W_i X_{ik} \right| \leq 2r \right).$$

$$\leq \mathbb{P}\left( \left| \frac{1}{n} \sum_{i=1}^{n} W_i X_{ik} \right| \geq \frac{r}{2} \right)$$

$$+ \mathbb{P}\left( \sum_{j=1}^{M} |\widetilde{\beta}_j - \beta_j^*| \left| \frac{1}{n} \sum_{i=1}^{n} X_{ij} X_{ik} \right| > \frac{r}{2} + (1 + \frac{1}{r})\epsilon \right)$$

$$+ \mathbb{P}\left( |S_n - B_n| \geq \frac{r}{2} + (1 + \frac{1}{r})\epsilon \right),$$

where the last inequality follows from the assumption that $\min_{j \in I^*} |\beta_j^*| \geq 3.5r + 3(1 + \frac{1}{r})\epsilon$. We bound the first term above using Hoeffding's inequality:

$$\mathbb{P}\left( \left| \frac{1}{n} \sum_{i=1}^{n} W_i X_{ik} \right| \geq r/2 \right) \leq \frac{\delta}{M}, \tag{4.15}$$

since, in particular, $r \geq 2L\sqrt{\frac{2\log(\frac{2M}{\delta})}{n}}$.

If *Condition Identif* holds, we can bound the second term of the last inequality of the display above by

$$\mathbb{P}\left( |\widetilde{\beta} - \beta^*|_1 > \frac{rk^*}{2d} + (1 + \frac{1}{r})\epsilon \right) \leq \frac{\delta}{M}, \tag{4.16}$$

as in Theorem 2.4, if $\frac{1}{2d} \geq \frac{4}{sb}$, with $b$ given by *Condition Stabil*. By Lemma 2.1, *Condition Identif* implies *Condition Stabil* with $b = 1 - d(7 + \epsilon)$, and so the restriction on $d$ is $d \leq \frac{s}{8+s(7+\epsilon)}$.

It remains to bound the term $\mathbb{P}(|S_n - B_n| \geq \frac{r}{2} + (1 + \frac{1}{r})\epsilon)$. For this, let $g(z) = e^z/(1 + e^z)$ and notice that Taylor's formula gives $g(u) - g(v) = g'(a)(u - v)$, for a point $a$ between $u$ and $v$, where $0 < g'(a) < 1$. Therefore

$$|S_n - B_n| \leq \sum_{j=1}^{M} |\widetilde{\beta}_j - \beta_j^*| \left| \frac{1}{n} \sum_{i=1}^{n} (1 - g'(a_i)) X_{ij} X_{ik} \right|, \qquad (4.17)$$

where $a_i$ is a point between $\sum_{j=1}^{M} \widetilde{\beta}_j X_{ij}$ and $\sum_{j=1}^{M} \beta_j^* X_{ij}$, for each $i$, and so

$$\left| a_i - \sum_{j=1}^{M} \beta_j^* X_{ij} \right| \leq L \sum_{j=1}^{M} |\widetilde{\beta}_j - \beta_j^*|, \quad \text{for each } i.$$

Let

$$G_n = \left\{ \sum_{j=1}^{M} |\widetilde{\beta}_j - \beta_j^*| \leq \frac{4rk^*}{sb} + (1 + \frac{1}{r})\epsilon) \right\}.$$

Therefore, by Theorem 2.4, for $b$ chosen as in the discussion following display (4.16) above, we have

$$\mathbb{P}(G_n^c) \leq \delta.$$

Notice that on the event $G_n$ we have

$$\left| a_i - \sum_{j=1}^{M} \beta_j^* X_{ij} \right| \leq \frac{4rLk^*}{sb} + L(1 + \frac{1}{r})\epsilon, \quad \text{for each } i.$$

This justifies the definition of the set $U$ in *Condition Lidentif*.

Combining the results above with (4.17) we obtain

$$\mathbb{P}\left( |S_n - B_n| \geq \frac{r}{2} + (1 + \frac{1}{r})d\epsilon \right)$$

$$\leq \mathbb{P}\left( \sum_{j=1}^{M} |\widetilde{\beta}_j - \beta_j^*| \left| \frac{1}{n} \sum_{i=1}^{n} (1 - g'(a_i)) X_{ij} X_{ik} \right| \geq \frac{r}{2} + (1 + \frac{1}{r})\epsilon \bigcap G_n \right)$$

$$+ \delta.$$

Note that if *Condition Identif* and *Lidentif* both hold for $d/2$ then

$$\left| \frac{1}{n} \sum_{i=1}^{n} (1 - g'(a_i)) X_{ij} X_{ik} \right| \leq \frac{d}{k^*}.$$

Thus, if $d \leq \frac{s}{16+2s(7+\epsilon)}$, and with $b$ chosen as in the discussion following display (4.16) above we have

$$\mathbb{P}\left( |S_n - B_n| \geq \frac{r}{2} + (1 + \frac{1}{r})d\epsilon \right) \leq \mathbb{P}(G_n^c \cap G_n) + \frac{\delta}{M} = \frac{\delta}{M}.$$

Therefore, collecting the bounds above, we obtain

$$\mathbb{P}(I^* \not\subseteq \hat{I}) \leq \frac{3k^*\delta}{M} \leq 3\delta.$$

The result for the $\ell_1 + \ell_2$ penalized estimator follows in an identical manner. By Lemma 4.3 in Appendix B below, if $\widetilde{\beta}_k = 0$ is a component of the solution $\widetilde{\beta}$ then

$$\left| \frac{1}{n} \sum_{i=1}^{n} X_{ik} \frac{\exp \sum_{j=1}^{M} \widetilde{\beta}_j X_{ij}}{1 + \exp \sum_{j=1}^{M} \widetilde{\beta}_j X_{ij}} - \frac{1}{n} \sum_{i=1}^{n} Y_i X_{ik} \right| \leq 2r.$$

Therefore the remainder of the proof is identical to the proof above, if we invoke Theorem 2.7 instead of Theorem 2.4. □

*Proof of Theorem 3.5.* In light of Proposition 3.3, it is enough to show that $\mathbb{P}(\hat{I} \subseteq I^*) \geq 1 - 2\delta$, for both the $\ell_1$ and $\ell_1 + \ell_2$ penalized least squares estimators.

We begin by showing that $\mathbb{P}(\hat{I} \subseteq I^*) \geq 2\delta$ for the $\ell_1$ penalized estimate. Let

$$h(\mu) = \frac{1}{n} \sum_{i=1}^{n} \left\{ Y_i - \sum_{j \in I^*} \mu_j X_{ij} \right\}^2 + 2r \sum_{j \in I^*} |\mu_j|,$$

and define

$$\hat{\mu} = \underset{\mu \in \mathbb{R}^{k^*}}{\arg\min} \, h(\mu). \tag{4.18}$$

Let

$$\mathcal{B} = \bigcap_{k \notin I^*} \left\{ \left| \frac{2}{n} \sum_{i=1}^{n} \left( Y_i - \sum_{j \in I^*} \hat{\mu}_j X_{ij} \right) X_{ik} \right| < 2r \right\}.$$

Let, by abuse of notation, $\hat{\mu} \in \mathbb{R}^M$ be the vector that has the components of $\hat{\mu}$ in positions corresponding to the index set $I^*$ and components equal to zero otherwise. By standard results in convex analysis, e.g. Lemma 4.1 in the Appendix B below it follows that, on the set $\mathcal{B}$, $\hat{\mu}$ is a solution of (2.2). Recall that $\hat{\beta}$ is a solution of (2.2) by construction. By definition $\hat{\beta}_k \neq 0$ for $k \in \hat{I}$. By construction, $\hat{\mu}_k \neq 0$ for $k \in S \subseteq I^*$, for some subset $S$. By Proposition 4.2 in Appendix B, any two solutions have non-zero elements in the same positions, therefore $\hat{I} = S \subseteq I^*$ on $\mathcal{B}$. Hence

$$\begin{aligned}
\mathbb{P}(\hat{I} \not\subseteq I^*) &\leq \mathbb{P}(\mathcal{B}^c) \leq \sum_{k \notin I^*} \mathbb{P} \left( \left| \frac{2}{n} \sum_{i=1}^{n} \left[ Y_i - \sum_{j \in I^*} \hat{\mu}_j X_{ij} \right] X_{ik} \right| \geq 2r \right) \\
&\leq \sum_{k \notin I^*} \mathbb{P} \left( \frac{1}{n} \left| \sum_{i=1}^{n} W_i X_{ik} \right| \geq \frac{r}{2} \right) \\
&+ \sum_{k \notin I^*} \mathbb{P} \left( \left| \sum_{j \in I^*} (\hat{\mu}_j - \beta_j^*) \left( \frac{1}{n} \sum_{i=1}^{n} X_{ij} X_{ik} \right) \right| \geq \frac{r}{2} \right)
\end{aligned}$$

$$\leq \sum_{k=1}^{M} \mathbb{P}\left(\frac{1}{n}\left|\sum_{i=1}^{n} W_i X_{ik}\right| \geq \frac{r}{2}\right)$$

$$+ \sum_{k \notin I^*} \mathbb{P}\left(\left|\sum_{j \in I^*}(\widehat{\mu}_j - \beta_j^*)\left(\frac{1}{n}\sum_{i=1}^{n} X_{ij}X_{ik}\right)\right| \geq \frac{r}{2}\right)$$

$$\leq \sum_{k=1}^{M} \mathbb{P}\left(\frac{1}{n}\left|\sum_{i=1}^{n} W_i X_{ik}\right| \geq \frac{r}{2}\right) + \sum_{k \notin I^*} \mathbb{P}\left(\sum_{j \in I^*}|\widehat{\mu}_j - \beta_j^*| \geq \frac{rk^*}{2d}\right),$$

where we used *Condition Identif* to obtain the last inequality. Recall now that if $Y \in \{0,1\}$ the choice

$$r \geq 2\sqrt{\frac{2\ln(\frac{2M}{\delta})}{n}},$$

guarantees, as in display (4.4) of the proof of Theorem 2.2, that

$$\sum_{k=1}^{M} \mathbb{P}\left(\frac{1}{n}\left|\sum_{i=1}^{n} W_i X_{ik}\right| \geq \frac{r}{2}\right) \leq \delta.$$

Repeating now the proof of Theorem 2.2, with $\widehat{\beta}$ replaced by $\widehat{\mu}$ and using only the variables corresponding to $I^*$, we obtain

$$|\widehat{\mu} - \beta^*|_1 \leq \frac{rk^*}{2d}$$

on the set

$$A_1 = \bigcap_{j \in I^*}\left\{\left|\frac{2}{n}\sum_{i=1}^{n} W_i X_{ij}\right| \leq r\right\}.$$

By Hoeffding's inequality

$$\mathbb{P}(A_1^c) \leq 2k^* \exp(-nr^2/8),$$

and the choice

$$r_{n,k^*}(\delta) \geq 2\sqrt{\frac{2\ln(\frac{2k^* M}{\delta})}{n}} \qquad (4.19)$$

implies that $\mathbb{P}(A_1^c) \leq \delta/M$, which in turn implies that

$$\sum_{k \notin I^*} \mathbb{P}\left(\sum_{j \in I^*}|\widehat{\mu}_j - \beta_j^*| \geq \frac{rk^*}{2d}\right) \leq \delta.$$

Here we used again the fact that, by Lemma 2.1, *Condition Identif* implies *Condition Stabil* and then reasoned as in Proposition 3.3 to conclude that the analogue of Theorem 2.2 can be used, for $d \leq \frac{1}{15}$.

The same conclusion holds if $Y \in \mathbb{R}$, by invoking Bernstein's inequality as in (4.5) and corresponding value of $r$ from the statement of this proposition, instead of Hoeffding's inequality.

Of course, the choice in (4.19) is not implementable, as $k^*$ is not known in practice, and we can always replaced it by a known upper bound, or the conservative bound $M$. This completes the proof for this part of the proposition.

It remains to show that $\mathbb{P}(\widehat{I} \subseteq I^*) \geq 1 - 2\delta$ for the $\ell_1 + \ell_2$ penalized estimate. We reason as above and let

$$m(\mu) = \frac{1}{n} \sum_{i=1}^{n} \{Y_i - \sum_{j \in I^*} \mu_j X_{ij}\}^2 + 2r \sum_{j \in I^*} |\mu_j| + c \sum_{j \in I^*} \mu_j^2,$$

and define

$$\widehat{\mu} = \operatorname*{arg\,min}_{\mu \in \mathbb{R}^{k^*}} m(\mu). \tag{4.20}$$

Then, by Lemma 4.3 in the Appendix B, $\widehat{b} = (\widehat{\mu}, 0)$, where 0 is a vector corresponding to indices in $I^{*c}$, is a solution of (2.3) on the set

$$\mathcal{B} = \bigcap_{k \in I^{*c}} \left\{ \left| \frac{2}{n} \sum_{i=1}^{n} \left( Y_i - \sum_{j \in I^*} \widehat{b}_j X_{ij} \right) X_{ik} \right| < 2r \right\}.$$

Recall that $\widehat{\beta}$ is a solution of (2.3) by construction, and that by Lemma 4.3 in the Appendix B, the solution is unique. Since, on the set $\mathcal{B}$, $\widehat{b}_k = 0$ for $k \in I^{*c}$, by construction, and $\widehat{\beta}_k = 0$ on $\widehat{I}^c$, by definition, we conclude that $\widehat{I} \subseteq I^*$ on the set $\mathcal{B}$. Therefore the proof is identical to the one above, where we now invoke *Condition Identif* with $d \leq \frac{1+c}{17.5}$ and the analogue of the proof of Theorem 2.3. □

*Proof of Theorem 3.7.* By Proposition 3.4, it is enough to show that $\mathbb{P}(\widehat{I} \subseteq I^*) \geq 1 - 2\delta$ for both the $\ell_1$ and $\ell_1 + \ell_2$ penalized estimate.

We begin by showing that $\mathbb{P}(\widehat{I} \subseteq I^*) \geq 1 - 2\delta$ for the $\ell_1$ penalized logistic regression estimate. Let

$$H(\mu) = \frac{1}{n} \sum_{i=1}^{n} \left\{ -Y_i \mu' \mathbf{X}_i + \log(1 + \exp \mu' \mathbf{X}_i) \right\} + 2r \sum_{j \in I^*} |\mu_j|,$$

and define

$$\widetilde{\mu} = \operatorname*{arg\,min}_{\mu \in \mathbb{R}^{k^*}} H(\mu). \tag{4.21}$$

Let

$$\mathcal{B}_1 = \bigcap_{k \notin I^*} \left\{ \left| \frac{1}{n} \sum_{i=1}^{n} X_{ik} \frac{\exp \sum_{j \in I^*} \widetilde{\mu}_j X_{ij}}{1 + \exp \sum_{j \in I^*} \widetilde{\mu}_j X_{ij}} - \frac{1}{n} \sum_{i=1}^{n} Y_i X_{ik} \right| \leq 2r \right\}.$$

Let, by abuse of notation, $\widetilde{\mu} \in \mathbb{R}^M$ be the vector that has the components of $\widetilde{\mu}$ in positions corresponding to the index set $I^*$ and components equal to zero otherwise. By standard results in convex analysis, e.g. Lemma 4.1 in the Appendix below it follows that, on the set $\mathcal{B}_1$, $\widetilde{\mu}$ is a solution of (2.2). Recall that $\widetilde{\beta}$ is a solution of (2.2) by construction. Then, by Proposition 4.2 in the Appendix B, any two solutions have non-zero elements in the same positions. Since, on the set $\mathcal{B}_1$, $\widehat{\beta}_k = 0$ for $k \in I^{*c}$ we conclude that $\widehat{I} \subseteq I^*$ on the set $\mathcal{B}_1$. Hence, reasoning as in Theorem 3.5 above

$$
\begin{aligned}
\mathbb{P}(\widehat{I} \nsubseteq I^*) \;\leq\; & \mathbb{P}(\mathcal{B}_1^c) \\
\leq\; & \sum_{k \notin I^*} \mathbb{P}\left( \frac{1}{n} \left| \sum_{i=1}^n W_i X_{ik} \right| \geq r \right) \\
& + \sum_{k \notin I^*} \mathbb{P}\left( \sum_{j \in I^*} |\widetilde{\mu}_j - \beta_j^*| \left| \frac{1}{n} \sum_{i=1}^n g'(a_i) X_{ij} X_{ik} \right| \geq r \right) \\
\leq\; & \sum_{k=1}^M \mathbb{P}\left( \frac{1}{n} \left| \sum_{i=1}^n W_i X_{ik} \right| \geq \frac{r}{2} \right) \\
& + \sum_{k \notin I^*} \mathbb{P}\left( \sum_{j \in I^*} |\widetilde{\mu}_j - \beta_j^*| \left| \frac{1}{n} \sum_{i=1}^n g'(a_i) X_{ij} X_{ik} \right| \geq \frac{r}{2} \right) \\
\leq\; & \delta + \sum_{k \notin I^*} \mathbb{P}\left( \sum_{j \in I^*} |\widetilde{\mu}_j - \beta_j^*| \geq r k^* / 2d \right) \\
\leq\; & \delta + \sum_{k \notin I^*} \mathbb{P}\left( \left\{ \sum_{j \in I^*} |\widetilde{\mu}_j - \beta_j^*| \geq r k^* / 2d \right\} \cap D_n \right) + \sum_{k \notin I^*} \mathbb{P}(D_n^c),
\end{aligned}
$$

where, as in (4.15), we used Hoeffding's inequality to bound by $\delta$ the first term, we used *Condition Lidentif* for the second term, and where

$$
D_n = \left\{ \sum_{j \in I^*} |\widetilde{\mu}_j - \beta_j^*| \leq \frac{4 r k^*}{sb} + \left(1 + \frac{1}{r}\right) \epsilon \right\},
$$

with $b = 1 - d(7 + \epsilon)$ and

$$
\epsilon = \frac{\log 2}{2^{(M \vee n)+1}} \times \frac{1}{r}.
$$

Notice that by the definition of $\epsilon$ and $r$, and since $0 < b, s \leq 1$, we always have $(1 + \frac{1}{r})\epsilon \leq \frac{r k^*}{sb}$. Thus, for our choice of $d$, we have $r k^*/2d > \frac{4 r k^*}{sb} + (1 + \frac{1}{r})\epsilon$. Therefore

$$
\mathbb{P}(\widehat{I} \nsubseteq I^*) \leq \delta + \sum_{k \notin I^*} \mathbb{P}(D_n \cap D_n^c) + \sum_{k \notin I^*} \mathbb{P}(D_n^c) = \delta + \sum_{k \notin I^*} \mathbb{P}(D_n^c) \leq \delta + \sum_{k=1}^M \mathbb{P}(D_n^c).
$$

Repeating now the proof of Theorem 2.4, with $\widetilde{\beta}$ replaced by $\widetilde{\mu}$, where we are now using only the variables corresponding to $I^*$ we obtain that

$$|\widetilde{\mu} - \beta^*|_1 \leq \frac{4}{b} r k^* + (1 + \frac{1}{r})\epsilon$$

on the set

$$A_2 = \left\{ \sup_{\beta \in \mathbb{R}^{k^*}} \frac{|(\mathbb{P}_n - \mathbb{P})(l(\beta^*) - l(\beta))|}{|\beta - \beta^*|_1 + \epsilon} \leq r \right\},$$

where, as in Theorem 2.4, we can show that $\mathbb{P}(A_2^c) \leq \delta/M$, for our choice of $r$. Therefore

$$\mathbb{P}(\hat{I} \not\subseteq I^*) \leq 2\delta.$$

It remains to show that the result above also holds for the $\ell_1 + \ell_2$ penalized estimator. Define

$$
\begin{aligned}
M(\mu) \;\; = \;\; & \frac{1}{n} \sum_{i=1}^{n} \left\{ -Y_i \sum_{j \in I^*} \mu_j X_{ij} + \log(1 + \exp \sum_{j \in I^*} \mu_j X_{ij}) \right\} \\
& + 2r \sum_{j \in I^*} |\mu_j| + c \sum_{j \in I^*} \mu_j^2,
\end{aligned}
\tag{4.22}
$$

and let

$$\widetilde{\mu} = \underset{\mu \in \mathbb{R}^{k^*}}{\arg \min} \, M(\mu). \tag{4.23}$$

Let

$$\mathcal{B}_1 = \bigcap_{k \notin I^*} \left\{ \left| \frac{1}{n} \sum_{i=1}^{n} X_{ik} \frac{\exp \sum_{j \in I^*} \widetilde{\mu}_j X_{ij}}{1 + \exp \sum_{j \in I^*} \widetilde{\mu}_j X_{ij}} - \frac{1}{n} \sum_{i=1}^{n} Y_i X_{ik} \right| \leq 2r \right\}.$$

Let, by abuse of notation, $\widetilde{\mu} \in \mathbb{R}^M$ be the vector that has the components of $\widetilde{\mu}$ in positions corresponding to the index set $I^*$ and components equal to zero otherwise. By standard Lemma 4.3 in the Appendix B below it follows that, on the set $\mathcal{B}_1$, $\widetilde{\mu}$ is a solution of (2.2). Recall that $\widetilde{\beta}$ is a solution of (2.2) by construction. Also by Lemma 4.3 in the Appendix B, the solution is unique. Since, on the set $\mathcal{B}_1$, $\widetilde{\beta}_k = 0$ for $k \in I^{*c}$ we conclude that $\hat{I} \subseteq I^*$ on the set $\mathcal{B}_1$. Therefore the remainder of the proof is identical to the one above. □

## Appendix B

### *4.1. Properties of $\ell_1$ penalized least squares and logistic regression solutions*

The solution of the $\ell_1$ penalized optimization problem may not be unique. However, in this case, all solutions have zero elements in the same positions, as we show below. We denote by $\mathcal{Y} = (Y_1, \ldots, Y_n)$ and by $\mathcal{X}$ the $n \times M$ matrix with entries $X_{ij}$. We let $L(\beta) = L(\mathcal{X}, \mathcal{Y}; \beta)$ be a function depending on the data and a parameter $\beta \in \mathbb{R}^M$. Let

$$\bar{\beta} = \arg\min_{\beta} L(\beta) + \lambda \sum_{j=1}^{M} |\beta_j| =: \arg\min_{\beta} f(\beta), \qquad (4.24)$$

for some fixed $\lambda > 0$. Let $S$ be the set of indices corresponding to the non-zero components of a solution $\bar{\beta}$ :

$$S = \{k : \bar{\beta}_k \neq 0, \ 1 \leq k \leq M\}.$$

**Lemma 4.1.** *If $L$ is differentiable in $\beta$ and if for any minima $\bar{\beta}^{(1)}, \bar{\beta}^{(2)}$*

$$\frac{\partial L(\bar{\beta}^{(1)})}{\partial \beta_j} = \frac{\partial L(\bar{\beta}^{(2)})}{\partial \beta_j}, \quad \textit{for all } 1 \leq j \leq M, \qquad (4.25)$$

*then all $\bar{\beta}$ satisfying (4.24) have non-zero components in the same positions.*

*Proof.* We recall that for any convex function $f : \mathbb{R}^M \to \mathbb{R}$ the subdifferential of $f$ at a point $\beta$ is the set $D_\beta = \{w \in \mathbb{R}^M : f(u) - f(\beta) \geq \langle w, u - \beta \rangle\}$. For the function $f$ defined in (4.24) this becomes

$$D_\beta = \{w \in \mathbb{R}^M : w = \nabla L(\beta) + \lambda v\},$$

where $\nabla L(\beta)$ is the $M$-dimensional vector having $\frac{\partial L(\beta)}{\partial \beta_j}$ as components and $v \in \mathbb{R}^M$ is such that

$$
\begin{array}{lll}
v_j & = & 1, & \text{if } \beta_j > 0 \\
v_j & = & -1, & \text{if } \beta_j < 0 \\
v_j & \in & [-1, 1], & \text{if } \beta_j = 0.
\end{array}
$$

By standard results in convex analysis, $\bar{\beta} \in \mathbb{R}^M$ is a point of local minimum for a convex function $f$ if and only if $0 \in D_{\bar{\beta}}$, where $0 \in \mathbb{R}^M$.

Therefore, $\bar{\beta}$ satisfies (4.24) if and only if

$$\left| \frac{\partial L(\bar{\beta})}{\partial \beta_j} \right| = \lambda |v|, \quad \text{for all} \quad 1 \leq j \leq M,$$

and so the index set $S$ of non-zero components of a solution is given by

$$S = \left\{ 1 \leq j \leq M : \left| \frac{\partial L(\bar{\beta})}{\partial \beta_j} \right| = \lambda \right\}.$$

Therefore, if (4.25) holds, $S$ is the same for all solutions. $\qquad \square$

**Proposition 4.2.** *Let $L$ correspond to either the least squares or the logistic criteria. Let $\bar{\beta}^{(1)}$ and $\bar{\beta}^{(2)}$ be two minima of (4.24). Then:*

*(1) $\mathcal{X}(\bar{\beta}^{(1)} - \bar{\beta}^{(2)}) = 0$, for either estimate.*

*(2) All solutions of (4.24), for either estimate, have non-zero components in the same positions.*

*Proof.* The proof uses simple properties of convex functions. First, we recall that the set of minima of a convex function is convex. Therefore, if $\bar{\beta}^{(1)}$ and $\bar{\beta}^{(2)}$ are two distinct points of minima, so is $\rho\bar{\beta}^{(1)} + (1-\rho)\bar{\beta}^{(2)}$, for any $0 < \rho < 1$. Re-write this convex combination as $\bar{\beta}^{2} + \rho\eta$, where $\eta = \bar{\beta}^{(1)} - \bar{\beta}^{(2)}$. Then, recall that the minimum value of any convex function is unique. For clarity, we argue separately for the two estimates.

$\ell_1$ *penalized least squares.* By the above arguments we have that

$$F(\rho) =: \frac{1}{n}\sum_{i=1}^{n}\left\{Y_i - (\bar{\beta}^{(2)} + \rho\eta)'X_i\right\}^2 + \lambda\sum_{j=1}^{M}|\bar{\beta}_j^{(2)} + \rho\eta_j| = c, \qquad (4.26)$$

where $c$ is some positive constant, for any $0 < \rho < 1$. By taking the derivative with respect to $\rho$ of $F(\rho)$ above we obtain

$$-\frac{2}{n}\sum_{i=1}^{n}Y_i\left(\sum_{j=1}^{M}\eta_j X_{ij}\right) + \frac{2}{n}\sum_{i=1}^{n}\left(\sum_{j=1}^{M}\eta_j X_{ij}\right)\left(\sum_{j=1}^{M}\beta_j X_{ij}\right)$$

$$+ \frac{2\rho}{n}\sum_{i=1}^{n}\left(\sum_{j=1}^{M}\eta_j X_{ij}\right)^2 + \lambda\sum_{j=1}^{M}\eta_j\mathrm{sign}(\bar{\beta}_j^{(2)} + \rho\eta_j) = 0.$$

Since the function $a + b\rho$ is continuous in $\rho$ then, on a small neighborhood $\mathcal{U}$ of $\rho$ the sign of $\bar{\beta}_j^{(2)} + \rho\eta_j$, for each $j$, will be constant. Therefore, on $\mathcal{U}$, the first two and the last term of the display of above are constant with respect to $\rho$. Denoting the sum of these terms by $C$ we have

$$\frac{2\rho}{n}\sum_{i=1}^{n}\left(\sum_{j=1}^{M}\eta_j X_{ij}\right)^2 + C = 0, \quad \text{for any} \quad \rho \in \mathcal{U}.$$

By taking again the derivative with respect to $\rho$ we obtain that $\mathcal{X}\eta = 0$, which is the result stated in the first part of this Lemma.

$\ell_1$ *penalized logistic regression.* We argue as above that the value of the function $G$ below evaluated at a point of minimum is constant, and we evaluate it at a convex combination of two minima, as before. Thus, defining $G(\rho)$ as the quantity below

$$\frac{1}{n}\sum_{i=1}^{n}\left\{-Y_i(\bar{\beta}^{(2)} + \rho\eta)'X_i + \log(1 + \exp(\bar{\beta}^{(2)} + \rho\eta)'X_i)\right\} + \lambda\sum_{j=1}^{M}|\bar{\beta}_j^{(2)} + \rho\eta_j|$$

we have that $G(\rho) = c$, for some positive constant $c > 0$. Reasoning as above, we can take the derivative of the above function twice, with respect to $\rho$. Then, on a small neighborhood $\rho \in \mathcal{V}$ we have

$$\frac{1}{n} \sum_{i=1}^{n} \left( \sum_{j=1}^{M} \eta_j X_{ij} \right)^2 \frac{\exp \sum_{j=1}^{M} (\beta_j + \rho \eta_j) X_{ij}}{1 + \exp \sum_{j=1}^{M} (\beta_j + \rho \eta_j) X_{ij}} = 0, \quad \text{for any} \;\; \rho \in \mathcal{V},$$

which implies that $\sum_{j=1}^{M} \eta_j X_{ij} = 0$ for all $i$, which in turn implies that $\mathcal{X}\eta = 0$, as claimed in part (1) of this proposition.

The second part of the proposition follows trivially from the first part and by Lemma 4.1. It is enough to show that $\mathcal{X}(\bar{\beta}^{(1)} - \bar{\beta}^{(2)}) = 0$ implies

$$\frac{\partial L(\bar{\beta}^{(1)})}{\partial \beta_j} = \frac{\partial L(\bar{\beta}^{(2)})}{\partial \beta_j}, \quad \text{for all } j.$$

For the $\ell_1$ penalized least squares estimator we have

$$\frac{\partial L(\beta)}{\partial \beta_j} = \frac{2}{n} \sum_{i=1}^{n} \left[ Y_i - \sum_{k=1}^{M} \beta_k X_{ik} \right] X_{ij} = \frac{2}{n} \sum_{i=1}^{n} Y_i X_{ij} - \frac{2}{n} \sum_{i=1}^{n} \sum_{k=1}^{M} \beta_k X_{ik} X_{ij},$$

and the last term is constant across solutions if $\mathcal{X}'(\bar{\beta}^{(1)} - \bar{\beta}^{(2)}) = 0$, for any two solutions, and this is implied by part (1).

For the $\ell_1$ penalized logistic regression estimate we have

$$\frac{\partial L(\beta)}{\partial \beta_k} = \frac{1}{n} \sum_{i=1}^{n} X_{ik} \frac{\exp \sum_{j=1}^{M} \beta_j X_{ij}}{1 + \exp \sum_{j=1}^{M} \beta_j X_{ij}} - \frac{1}{n} \sum_{i=1}^{n} Y_i X_{ik}.$$

This will be constant across solutions if $\sum_{j=1}^{M} \beta_j X_{ij}$, for all $i$, is the same for all solutions, which is again implied by the result in part (1). This concludes the proof of this proposition. $\qquad \square$

## 4.2. Properties of the $\ell_1 + \ell_2$ penalized least squares and logistic regression solutions

We discuss below a number of properties of the solution of the $\ell_1 + \ell_2$ penalized optimization problem. We begin by giving this result in terms of general likelihood functions and we obtain the results for our two examples as consequences. As in the previous sub-section, we let $L(\beta) = L(\mathcal{X}, \mathcal{Y}; \beta)$ be any function depending on the data and a parameter $\beta \in \mathbb{R}^M$. Let

$$\bar{\beta} = \arg \min_{\beta} L(\beta) + \lambda \sum_{j=1}^{M} |\beta_j| + c \sum_{j=1}^{M} \beta_j^2 =: \arg \min_{\beta} s(\beta), \qquad (4.27)$$

for some given tuning parameters $\lambda, c > 0$. We note that this solution is different than the one introduced in the previous subsection, but for brevity we use the same notation.

**Lemma 4.3.** *If L is differentiable in β then a solution of (4.27) satisfies*

$$\left| \frac{\partial L(\bar{\beta})}{\partial \beta_j} + 2c\bar{\beta}_j \right| = \lambda, \quad if \bar{\beta}_j \neq 0, \tag{4.28}$$

$$\left| \frac{\partial L(\bar{\beta})}{\partial \beta_j} + 2\gamma\bar{\beta}_j \right| = \left| \frac{\partial L(\bar{\beta})}{\partial \beta_j} \right| \leq \lambda, \quad if \bar{\beta}_j = 0.$$

*Moreover, the solution of (4.27) is unique for both the square and the logistic losses, respectively.*

*Proof.* Appealing to the elementary properties of convex functions introduced in Lemma 4.1 and applying them now to the function $s$ above we trivially obtain the first part of this Lemma.

For the moreover part, let $\bar{\beta}^{(1)}$ and $\bar{\beta}^{(2)}$ be two solutions of (4.27). We show below that $\bar{\beta}^{(1)} = \bar{\beta}^{(2)}$ for the two losses under study. Since $s$ is a convex function of $\beta$, for either loss, any convex combination of solutions is a solution, and $s(\beta)$ is constant across solutions. Consider as before the convex combination $\bar{\beta}^2 + \rho\eta$, where $\eta = \bar{\beta}^{(1)} - \bar{\beta}^{(2)}$. Recall that the minimum value of any convex function is unique. Then, for the $\ell_1 + \ell_2$ penalized least square estimator we obtain:

$$F_1(\rho) =: \frac{1}{n} \sum_{i=1}^{n} \left\{ Y_i - (\bar{\beta}^{(2)} + \rho\eta)' X_i \right\}^2 + \lambda \sum_{j=1}^{M} |\bar{\beta}_j^{(2)} + \rho\eta_j| + \gamma \sum_{j=1}^{M} (\bar{\beta}_j^{(2)} + \rho\eta_j)^2 = c,$$

where $c$ is some positive constant, for any $0 < \rho < 1$. Reasoning now exactly as in Proposition 4.2 above and taking the derivative with respect to $\rho$ twice, we obtain

$$\frac{2}{n} \sum_{i=1}^{n} \sum_{j=1}^{M} \eta_j X_{ij})^2 + 2\gamma \sum_{j=1}^{M} \eta_j^2 = 0,$$

which immediately implies $\eta_j = 0$ for all $j$, that is $\bar{\beta}^{(1)} = \bar{\beta}^{(2)}$.

The same conclusion can be obtained for the logistic regression estimator, where we now differentiate twice with respect to $\rho$ the function $G_1(\rho) = G(\rho) + \gamma \sum_{j=1}^{M} (\bar{\beta}_j^{(2)} + \rho\eta_j)^2$, with $G(\rho)$ defined in display (4.27) of Proposition 4.2. This yields

$$\frac{1}{n} \sum_{i=1}^{n} \left( \sum_{j=1}^{M} \eta_j X_{ij} \right)^2 \frac{\exp \sum_{j=1}^{M} (\beta_j + \rho\eta_j) X_{ij}}{1 + \exp \sum_{j=1}^{M} (\beta_j + \rho\eta_j) X_{ij}} + 2\gamma \sum_{j=1}^{M} \eta_j^2 = 0.$$

Reasoning as above we again obtain $\bar{\beta}^{(1)} = \bar{\beta}^{(2)}$. This completes the proof of this Lemma. □

## Acknowledgements

I am grateful to Ingo Ruczinski, Sara van de Geer, Sasha Tsybakov, Vladimir Koltchinskii and Adrian Barbu for inspiring conversations.